\documentclass[reqno,11pt]{amsart}
\usepackage[margin=1.15in]{geometry}
\usepackage{amsmath,amssymb,mathtools}
\usepackage{enumitem,cite}
\usepackage{tikz}
\usepackage{hyperref}
\hypersetup{pdftitle={Nonradial entire solutions of the critical Henon equation for every alpha greater than two},pdfauthor={Phuong Le},pdflang={en-US}}
\allowdisplaybreaks
\newtheorem{theorem}{Theorem}[section]
\newtheorem{proposition}[theorem]{Proposition}
\newtheorem{lemma}[theorem]{Lemma}
\newtheorem{corollary}[theorem]{Corollary}
\theoremstyle{definition}
\newtheorem{definition}[theorem]{Definition}
\theoremstyle{remark}
\newtheorem{remark}[theorem]{Remark}
\title[Nonradial solutions of the critical H\'enon equation]{Nonradial entire solutions of the critical H\'enon equation for every \(\alpha>2\)}
\author{Phuong Le}
\address{Faculty of Economic Mathematics, University of Economics and Law, Ho Chi Minh City, Vietnam; Vietnam National University, Ho Chi Minh City, Vietnam}
\email{phuongl@uel.edu.vn}
\thanks{ORCID: \href{https://orcid.org/0000-0003-4724-7118}{0000-0003-4724-7118}.}
\begin{document}
\begin{abstract}
We consider the critical H\'enon equation \(-\Delta u=|x|^{\alpha}u^{(N+2+2\alpha)/(N-2)}\) in \(\mathbb{R}^N\), \(N\ge3\). Gladiali, Grossi, and Neves (\emph{Advances in Mathematics}, 2013) constructed nonradial continua bifurcating at \(\alpha_k=2(k-1)\), \(k\ge2\), and conjectured that these continua are confined to the even parameters. Jiang and Xiong (2026) disproved this conjecture by obtaining nonradial solutions in punctured neighborhoods of \(\alpha_k\) for even \(k>(N-2)/2\). We prove that positive nonradial finite-energy entire solutions exist for every \(\alpha>2\) in every dimension \(N\ge3\). The proof introduces a global variational construction by direct minimization on a monotonicity cone on the Emden--Fowler cylinder. Angular and axial order restore compactness in every finite \(L^r\), \(r>2\), despite Sobolev supercriticality, and invariance of the linear resolvent yields the full equation. This constructs a solution at each prescribed parameter without continuing local bifurcation branches. The solutions have Newtonian decay, are invariant under \(O(N-1)\times O(1)\) and the Kelvin transform, and have energy strictly below the radial level. We also obtain a second solution for every \(\alpha>6\) in even dimension, nonvertical branches at every \(\alpha_k\) with \(k\ge2\) even without dimension restrictions, a cubic energy law, and rigidity near \(\alpha=0\) under bounded cylindrical amplitude.
\end{abstract}
\maketitle

\section{Introduction and main results}

We consider positive solutions of
\begin{equation}\label{eq:henon}
 -\Delta u=|x|^\alpha u^{p_\alpha}\quad\hbox{in }\mathbb{R}^N,
 \qquad p_\alpha=\frac{N+2+2\alpha}{N-2},\qquad N\ge3.
\end{equation}
The weight originates in H\'enon's model of spherical stellar systems \cite{He}.
For the Dirichlet problem on a ball, Ni \cite{Ni} showed that it permits positive radial solutions beyond the ordinary Sobolev exponent. Symmetry breaking for ground states of the subcritical ball problem was studied by Smets, Willem, and Su \cite{SWS}.
Here we work on the whole space at the weighted critical exponent, which balances the radial weight under the dilation
\(u(x)\mapsto \lambda^{(N-2)/2}u(\lambda x)\).
The positive radial solutions in \(D^{1,2}(\mathbb{R}^N)\) are dilations of
\begin{equation}
 U_\alpha(x)=
 [(N+\alpha)(N-2)]^{1/(p_\alpha-1)}
 (1+|x|^{\alpha+2})^{-(N-2)/(\alpha+2)}.
\end{equation}
When the weight is singular, the radial solutions are the only ones: for \(-2<\alpha<0\), Chou and Chu \cite{CC} proved that a positive weak solution in \(D^{1,2}(\mathbb R^N)\), the equation being required to hold weakly also at the origin, is a dilation \(\lambda^{(N-2)/2}U_\alpha(\lambda x)\) of \(U_\alpha\). At \(\alpha=0\), translation invariance allows the center of symmetry to vary, and the positive entire solutions are the translates and dilates of a single bubble \cite{CGS}. When \(\alpha>0\), the exponent \(p_\alpha\) lies above the Sobolev exponent, and the energy functional of \eqref{eq:henon} is not even defined on all of \(D^{1,2}(\mathbb{R}^N)\).

Gladiali, Grossi, and Neves \cite{GGN} identified the additional angular degeneracy at \(\alpha_k=2(k-1)\), \(k\ge2\), and constructed continua of positive nonradial solutions bifurcating there from the radial family. They conjectured that the nonradial solutions in these continua occur only at even values of \(\alpha\); see the discussion following \cite[Theorem~1.6]{GGN}.
The corresponding bifurcation problem for critical \(p\)-Laplacian H\'enon equations was studied by Dai, Duan, Gui, and Li \cite{DDGL}; when \(p=2\), their distinguished parameters reduce to the same values \(2(k-1)\).

Jiang and Xiong \cite{JX} recently disproved this conjecture. For every even \(k>(N-2)/2\), they computed a nonzero bifurcation slope and obtained positive nonradial solutions for every \(\alpha\) in a punctured neighborhood of \(\alpha_k\). Their paper also gives a global continuum alternative, but it does not establish existence at every prescribed \(\alpha>2\).

Our main result establishes existence on the full half-line \(\alpha>2\), in every dimension \(N\ge3\), with energy strictly below the radial level. A global variational argument constructs a solution directly at each prescribed parameter, independently of the continuation or eventual shape of the local branches.

The unrestricted quotient already exhibits the difficulty. Fix \(x_0\ne0\) and a nonzero nonnegative \(\chi\in C_c^\infty(B_1)\), and set
\(u_\varepsilon(x)=\varepsilon^{-(N-2)/2}\chi((x-x_0)/\varepsilon)\).
Its Dirichlet integral is constant, while
\[
 \int_{\mathbb{R}^N}|x|^\alpha u_\varepsilon^{p_\alpha+1}\,dx
 =\varepsilon^{-\alpha}
 \left(|x_0|^\alpha\int_{\mathbb{R}^N}\chi^{p_\alpha+1}+o(1)\right).
\]
Consequently
\begin{equation}
 \frac{\int_{\mathbb{R}^N}|\nabla u_\varepsilon|^2}
 {\left(\int_{\mathbb{R}^N}|x|^\alpha u_\varepsilon^{p_\alpha+1}\right)^{2/(p_\alpha+1)}}
 \sim C\varepsilon^{(N-2)\alpha/(N+\alpha)}
 \longrightarrow0.
\end{equation}
Thus the unrestricted quotient has infimum zero, and direct minimization is impossible. Axial symmetry alone does not help, since concentration can still occur on the axis, whose \(O(N-1)\)-orbits are points.

To make direct minimization possible, we impose both angular and axial order on the cylindrical profiles.
Monotonicity constraints for supercritical Neumann problems were developed by Serra--Tilli \cite{ST} and Bonheure--Noris--Weth \cite{BNW}. Angular monotonicity cones for supercritical annular problems appear in \cite{BCNW,CM}. Bonheure, Casteras, and Gladiali \cite{BCG} used cones of functions monotone in the polar angle, together with Kelvin invariance, to separate global branches of nonradial solutions bifurcating from the radial one for a Hardy--Sobolev equation in \(\mathbb R^N\). Here the cone serves a different purpose: we minimize on it directly at each fixed parameter, and an additional monotonicity in the axial variable of the cylinder restores the compactness that the angular order alone does not provide. Aghajani and Cowan \cite{AC} study a supercritical H\'enon equation on a cylinder.
For the present problem, angular order toward the equator prevents concentration on the axis. Away from that axis, an axially symmetric profile depends on two variables, so every finite exponent is available.
The Emden--Fowler cylinder is unbounded, and we impose an additional evenness and monotonicity in its axial variable. A quantitative tail estimate then gives strong compactness in every finite \(L^r\), \(r>2\).
To recover the full Euler equation, we express the two orders through translation and rotation Killing fields. Commutation with the linear inverse proves cone invariance at the \(L^2\) level, and a Hilbert-space identity removes the constraint.

Set
\[
 m=\frac{N-2}{2},\qquad \mathcal{C}=\mathbb{R}\times\mathbb{S}^{N-1},
 \qquad u(e^t\theta)=e^{-mt}w(t,\theta).
\]
We write \(\vartheta=\arccos\theta_N\) for the colatitude from the north pole. Reflection across the equator restricts the angular coordinate to \(0\le\vartheta\le\pi/2\).

\begin{theorem}[Existence on the full half-line]\label{thm:henon}
Let \(N\ge3\) and \(\alpha>2\). Equation \eqref{eq:henon} has a positive nonradial solution
\[
 u\in C^2(\mathbb{R}^N)\cap D^{1,2}(\mathbb{R}^N).
\]
It is \(O(N-1)\times O(1)\)-invariant and satisfies
\[
 u(x)=|x|^{2-N}u\!\left(\frac{x}{|x|^2}\right)\quad(x\ne0).
\]
Its cylindrical profile is even in \(t\), with
\[
 w_t<0\quad(t>0),\qquad
 w_\vartheta>0\quad(t\in\mathbb{R},\ 0<\vartheta<\pi/2).
\]
For some \(A>0\),
\[
 \lim_{|x|\to0}u(x)=A,\qquad
 \lim_{|x|\to\infty}|x|^{N-2}u(x)=A
\]
uniformly in direction. Moreover,
\begin{equation}\label{eq:energy-identity}
 \int_{\mathbb{R}^N}|\nabla u|^2\,dx
 =\int_{\mathbb{R}^N}|x|^\alpha u^{p_\alpha+1}\,dx<\infty.
\end{equation}
\end{theorem}

The construction gives more information than nonradiality alone. The next result identifies the first angular term and compares the solution with the radial family. For functions with both integrals finite, define
\[
 \mathcal E_\alpha(v)
 =\frac12\int_{\mathbb{R}^N}|\nabla v|^2
 -\frac1{p_\alpha+1}\int_{\mathbb{R}^N}|x|^\alpha v^{p_\alpha+1}.
\]

\begin{theorem}[Energy and the first anisotropic term]\label{thm:geometry}
The solution in Theorem~\ref{thm:henon} can be chosen so that
\[
 \mathcal E_\alpha(u)<\mathcal E_\alpha(U_\alpha).
\]
There is \(B>0\) such that
\begin{align}
 u(x)&=A+B\bigl(|x|^2-Nx_N^2\bigr)+O(|x|^4)
       &&(|x|\to0),\label{eq:origin-expansion}\\
 u(r\theta)&=Ar^{2-N}
      +Br^{-N}(1-N\theta_N^2)+O(r^{-N-2})
       &&(r\to\infty),\label{eq:infinity-expansion}
\end{align}
uniformly in direction. In particular,
\[
 D^2u(0)=2B\,\operatorname{diag}(1,\ldots,1,-(N-1)),
\]
so the origin is a nondegenerate saddle point. At infinity,
\[
 \nabla u(x)=-(N-2)A\frac{x}{|x|^N}+O(|x|^{-N-1}).
\]
\end{theorem}

The energy gap and the coefficients \(A,B\) remain bounded away from zero, uniformly over the cone-minimizing solutions when \(\alpha\) ranges over a compact subinterval of \((2,\infty)\). This follows from compactness under parameter variation; see Corollary~\ref{cor:uniform-parameters}.

The same construction works for other symmetry groups, as long as the profile can be required to be largest where the orbits have dimension \(N-2\). In even dimension there is a second such group, and it gives a second solution once \(\alpha>6\), the value at which the radial solution acquires a negative direction of degree four \cite{GGN}. Let \(N=2d\) with \(d\ge2\), write \(x=(x',x'')\) with \(x',x''\in\mathbb{R}^d\), and let \(G_d\) be the group generated by \(O(d)\times O(d)\) and the exchange \((x',x'')\mapsto(x'',x')\). The polynomial
\[
 H_4(x)=|x'|^4+|x''|^4-\frac{2(d+2)}{d}|x'|^2|x''|^2
\]
is harmonic and \(G_d\)-invariant.

\begin{theorem}[A second solution in even dimension]\label{thm:second}
Let \(N=2d\ge4\) and \(\alpha>6\). Equation \eqref{eq:henon} has a positive \(G_d\)-invariant solution \(u_2\in C^2(\mathbb{R}^N)\cap D^{1,2}(\mathbb{R}^N)\) with the following properties.
\begin{enumerate}[label=\textup{(\roman*)}]
\item \(u_2(x)=|x|^{2-N}u_2(x/|x|^2)\) for \(x\ne0\), and \(\mathcal E_\alpha(u_2)<\mathcal E_\alpha(U_\alpha)\).
\item For each \(r>0\), the restriction of \(u_2\) to \(\{|x|=r,\ |x''|<|x'|\}\) is a strictly increasing function of \(|x''|/|x'|\).
\item There are \(A_2,B_2>0\) such that, uniformly in direction,
\begin{align*}
 u_2(x)&=A_2-B_2H_4(x)+O(|x|^8)&&(|x|\to0),\\
 u_2(r\theta)&=A_2r^{2-N}-B_2r^{-N-2}H_4(\theta)+O(r^{-N-6})&&(r\to\infty).
\end{align*}
\end{enumerate}
In particular \(D^2u_2(0)=0\), so \(u_2\) is not obtained from the solution of Theorem~\ref{thm:henon} by a rotation or a dilation.
\end{theorem}

The threshold \(\alpha=2\) in the first construction comes from an admissible degree-two variation of the radial profile. We next cross that threshold locally. Put
\[
 W_\alpha(t)=e^{mt}U_\alpha(e^t).
\]
For the local argument we work with relative perturbations \(w=W_\alpha(1+h)\), \(\|h\|_\infty<1/2\). This keeps the profile positive and makes the nonlinear map smooth for every \(p_\alpha>1\). The condition \(k>(N-2)/2\) in \cite{JX} is equivalent to \(p_{\alpha_k}>3\), where \(\alpha_k=2(k-1)\), and the relative formulation does not need it. For \(k=2\), for instance, the theorem of Jiang and Xiong applies when \(N=3,4,5\).

\begin{theorem}[Branches at the even degenerate values]\label{thm:local-two}
Let \(N\ge3\) and let \(k\ge2\) be even. There is \(\varepsilon_{N,k}>0\) such that, for every \(0<|\alpha-\alpha_k|<\varepsilon_{N,k}\), equation \eqref{eq:henon} has a positive nonradial \(O(N-1)\times O(1)\)-invariant solution in \(C^2(\mathbb{R}^N)\cap D^{1,2}(\mathbb{R}^N)\), with Newtonian decay and the Kelvin identity of Theorem~\ref{thm:henon}. These solutions lie on a smooth local cylindrical branch through \(W_{\alpha_k}\). On the part of the branch with \(\alpha<\alpha_k\) their energy is strictly greater than that of \(U_\alpha\), and on the part with \(\alpha>\alpha_k\) it is strictly smaller.
\end{theorem}

For \(k=2\) this covers an interval below \(2\) in every dimension, and for every even \(k\) it removes the condition \(k>(N-2)/2\) of \cite{JX}. Corollary~\ref{cor:branch-structure} gives the local uniqueness, an explicit cubic energy law, and the Morse index in the even axial class. The same cubic coefficient gives a lower bound for the energy gap of the global cone minimizers near \(2\), without assuming that they approach the radial profile; see Corollary~\ref{cor:cubic-gap}.

The existence results proved here do not cover the whole interval \((0,2)\), and unrestricted existence at \(\alpha=2\) in odd dimension remains unresolved. A companion study \cite{LeConc} develops concentration constructions near zero and near selected interior parameters. Those constructions are independent of the variational and local arguments in this paper. The following statements describe the endpoint rigidity relevant to the present construction.

\begin{proposition}[Local rigidity]\label{prop:local-rigidity}
Consider positive solutions with an even cylindrical profile.
\begin{enumerate}[label=(\roman*)]
\item If \(N\ge3\) is odd, there is \(\rho_N>0\) such that, at \(\alpha=2\),
\[
 \left\|\frac{w}{W_2}-1\right\|_\infty<\rho_N
 \quad\Longrightarrow\quad w=W_2.
\]
\item For every \(N\ge3\), there are \(\delta_N,\rho_N>0\) such that, for \(0<\alpha<\delta_N\),
\[
 \left\|\frac{w}{W_\alpha}-1\right\|_\infty<\rho_N
 \quad\Longrightarrow\quad w=W_\alpha.
\]
\end{enumerate}
No angular symmetry is assumed in these two assertions.
\end{proposition}

An even profile fixes the dilation phase, and Lemma~\ref{lem:t-symmetry} shows that every positive entire solution with Newtonian decay can be brought into this normalization by a dilation. Proposition~\ref{prop:local-rigidity} is a statement about a neighborhood of the radial solution in the relative norm, not a nonexistence theorem.

The conclusion near zero extends to every fixed bound on the cylindrical amplitude. For a positive entire solution with Newtonian decay, put
\[
 \mathcal A(u)=\sup_{x\ne0}|x|^m u(x)=\|w\|_{L^\infty(\mathcal C)}.
\]
This quantity is unchanged by dilation.

\begin{theorem}[Bounded-amplitude rigidity near zero]\label{thm:amplitude}
Let \(N\ge3\). For every \(M>0\), there is \(\delta_{N,M}>0\) such that every positive classical entire solution of \eqref{eq:henon} with Newtonian decay, \(0<\alpha<\delta_{N,M}\), and \(\mathcal A(u)\le M\) is radial. In particular, if \(\alpha_j\downarrow0\) and \(u_j\) are positive nonradial entire solutions with Newtonian decay, then
\begin{equation}\label{eq:amplitude-blowup}
 \sup_{x\ne0}|x|^m u_j(x)\longrightarrow\infty.
\end{equation}
No angular symmetry or energy bound is assumed.
\end{theorem}

The proof uses conservation of the cylindrical Hamiltonian to rule out a nonzero stationary tail in a bounded-amplitude limit. This gives uniform decay and relative compactness after dilation. Combining it with the balance identity and local rigidity yields Theorem~\ref{thm:amplitude}. In particular, the compact cone minimizers are radial for all sufficiently small positive \(\alpha\); see Corollary~\ref{cor:small-cone}.

The two constructions get compactness in different ways. The global one uses the angular order toward the equator and the axial order toward the center of the cylinder. The local one uses the exponential decay of the sources and works in a relative neighborhood of the radial profile. Its solutions with \(\alpha<2\) are not in the cone: they increase toward the poles.

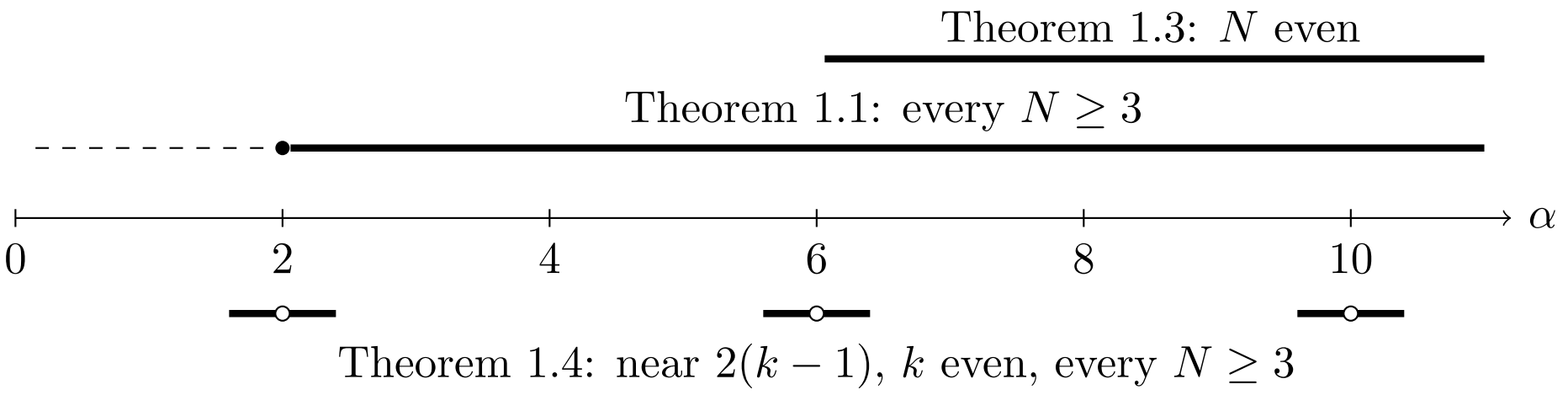
\begin{figure}[ht]
\centering
\begin{tikzpicture}[x=1.05cm,y=1cm,font=\small]
\draw[->] (0,0) -- (11.2,0) node[right] {\(\alpha\)};
\foreach \x in {0,2,4,6,8,10} {\draw (\x,0.07) -- (\x,-0.07) node[below] {\(\x\)};}
\draw[line width=1.6pt] (2.06,0.55) -- (11,0.55);
\node[above] at (6.5,0.55) {Theorem~\ref{thm:henon}: every \(N\ge3\)};
\draw[line width=1.6pt] (6.06,1.25) -- (11,1.25);
\node[above] at (8.5,1.25) {Theorem~\ref{thm:second}: \(N\) even};
\foreach \x in {2,6,10} {\draw[line width=1.6pt] (\x-0.4,-0.75) -- (\x+0.4,-0.75); \fill[white,draw=black] (\x,-0.75) circle (1.6pt);}
\node[below] at (6,-0.85) {Theorem~\ref{thm:local-two}: near \(2(k-1)\), \(k\) even, every \(N\ge3\)};
\draw[dashed] (0.15,0.55) -- (1.94,0.55);
\fill (2,0.55) circle (1.6pt);
\end{tikzpicture}
\caption{Nonradial solutions of \eqref{eq:henon}. Solid segments: existence. The intervals around \(2(k-1)\) have lengths depending on \(N\) and \(k\). At the dot \(\alpha=2\) there is an explicit family in even dimension (Section~\ref{sec:local-endpoints}) and local rigidity in odd dimension (Proposition~\ref{prop:local-rigidity}). The dashed segment indicates parameters not covered by the existence results of this paper; the companion study treats further concentration regimes. As \(\alpha\downarrow0\), the cylindrical amplitude of every nonradial family diverges (Theorem~\ref{thm:amplitude}).}
\label{fig:alpha}
\end{figure}

Figure~\ref{fig:alpha} summarizes the existence statements proved in this paper.

The proof is organized around a cylindrical equation with an arbitrary positive mass. Sections~\ref{sec:compactness}--\ref{sec:henon-proof} establish the cone construction, its energy gap, and its end expansions, and Section~\ref{sec:second} adapts it to \(G_d\). Section~\ref{sec:parameters} gives compactness under parameter variation. Sections~\ref{sec:relative}--\ref{sec:local-branch} construct the branches, and Section~\ref{sec:local-endpoints} treats local rigidity. The last two sections contain the application to a Hardy potential and the obstructions for \(\alpha\le2\).

\section{A cylindrical equation and its monotonicity cone}

For \(\mu>0\), put
\[
 P_\mu=-\partial_t^2-\Delta_{\mathbb{S}^{N-1}}+\mu^2,\qquad
 \|w\|_\mu^2=\int_\mathcal{C}
 \bigl(|w_t|^2+|\nabla_{\mathbb{S}^{N-1}}w|^2+\mu^2w^2\bigr).
\]
This norm is equivalent to the usual \(H^1(\mathcal{C})\)-norm.
The angular measure for equator-even axially symmetric functions is
\(2|\mathbb{S}^{N-2}|\sin^{N-2}\vartheta\,d\vartheta\) on
\(0<\vartheta<L:=\pi/2\). Also,
\begin{equation}
 \Delta_{\mathbb{S}^{N-1}}w
 =w_{\vartheta\vartheta}+(N-2)\cot\vartheta\,w_\vartheta.
\end{equation}

\begin{definition}\label{def:cone}
The cone \(\mathcal{K}\subset H^1(\mathcal{C})\) consists of nonnegative
\(O(N-1)\times O(1)\)-invariant functions satisfying
\[
 w(-t,\vartheta)=w(t,\vartheta),\qquad
 w_t\le0\ \ (t>0),\qquad
 w_\vartheta\ge0\ \ (0<\vartheta<L).
\]
The derivative inequalities are understood distributionally on interior coordinate strips.
\end{definition}

The cone is closed and convex. The symmetry conditions are closed linear conditions, while positivity and the derivative inequalities pass to an \(H^1\) limit by testing on compact strips. Consequently \(\mathcal{K}\) is weakly closed. It contains positive radial profiles that are even and decreasing in \(|t|\).

The angular order has a coordinate-free expression that will be useful for the inverse.
Let
\[
 Z=\theta_1\partial_{\theta_N}-\theta_N\partial_{\theta_1},
 \qquad Q=\{\theta\in\mathbb{S}^{N-1}:\theta_1>0,\ \theta_N>0\}.
\]
This is the rotation Killing field in the \((\theta_1,\theta_N)\)-plane.
Both \(Z\) and \(\partial_t\) are divergence free, have length at most one, and commute with \(P_\mu\).
For an axially symmetric profile, write \(w(t,\theta)=V(t,s)\), \(s=\theta_N\). Then
\begin{equation}\label{eq:killing-order}
 Zw=\theta_1V_s=-\frac{\theta_1}{\sin\vartheta}w_\vartheta.
\end{equation}
Thus \(w_\vartheta\ge0\) on \(0<\vartheta<L\) is equivalent to \(Zw\le0\) on \(\mathbb{R}\times Q\), also at the distributional level.

The following cylindrical result is the main construction used in the paper.

\begin{theorem}[Cone minimizers on the cylinder]\label{thm:cylinder}
Let \(N\ge3\), \(\mu>0\), and \(p>1\). There is a positive solution
\(w\in\mathcal{K}\cap C^\infty(\mathcal{C})\) of
\begin{equation}\label{eq:cylinder}
 P_\mu w=w^p.
\end{equation}
It minimizes the quotient
\begin{equation}\label{eq:quotient}
 S(\mu,p)=
 \inf_{\substack{v\in\mathcal{K}\\v\ne0}}
 \frac{\|v\|_\mu^2}
 {\left(\int_\mathcal{C} v^{p+1}\right)^{2/(p+1)}}.
\end{equation}
If
\begin{equation}\label{eq:breaking-condition}
 \frac{\mu^2(p+1)^2}{4}>\mu^2+2N,
\end{equation}
then every such minimizing solution is nonradial, has strictly positive angular derivative on \(0<\vartheta<L\), and has strictly negative axial derivative for \(t>0\).
\end{theorem}

The theorem is proved in Sections~\ref{sec:compactness}--\ref{sec:variation}.
The cone controls two distinct defects of compactness: concentration at the angular poles and translation along the cylinder. Its role in the Euler equation requires a separate argument, since a cone minimizer need not be a critical point of a functional on the full Hilbert space.

\section{Integrability and compactness}\label{sec:compactness}

On any angular strip separated from the poles, an axially symmetric profile is a function of two variables. Angular monotonicity controls the omitted polar caps by values on such a strip. This gives integrability at every finite exponent.

\begin{lemma}[All finite exponents]\label{lem:embedding}
For every \(2\le r<\infty\), there is \(C_{\mu,r}>0\) such that
\[
 \|w\|_{L^r(\mathcal{C})}\le C_{\mu,r}\|w\|_\mu
 \quad(w\in\mathcal{K}).
\]
For \(r>2\), every \(H^1\)-bounded sequence in \(\mathcal{K}\) has a subsequence converging strongly in \(L^r(\mathcal{C})\).
\end{lemma}

\begin{proof}
Set \(a=\pi/4\) and \(D=\mathbb{R}\times(a,L)\). The angular weight is bounded above and below on \(D\), so
\[
 \|w\|_{H^1(D)}\le C_\mu\|w\|_\mu.
\]
For almost every \(t\), monotonicity gives
\[
 w(t,\vartheta)^r
 \le\frac1{L-a}\int_a^Lw(t,y)^r\,dy
 \quad(0<\vartheta<a).
\]
After integration,
\[
 \int_\mathcal{C} w^r\le C\int_Dw^r.
\]
The ordinary two-dimensional embedding \(H^1(D)\hookrightarrow L^r(D)\), valid for every finite \(r\ge2\), proves the norm estimate.

For \(0<\delta<a\), the same argument yields the polar-cap bound
\begin{equation}\label{eq:cap}
 \int_\mathbb{R}\int_0^\delta
 w^r\sin^{N-2}\vartheta\,d\vartheta\,dt
 \le C_{\mu,r}\delta^{N-1}\|w\|_\mu^r.
\end{equation}
The south-pole cap is identical by reflection.

To control the axial ends, let \(b(t)=w(t,L)\) be the angular trace.
The one-dimensional trace estimate in the angular variable, integrated in \(t\), gives
\begin{equation}
 \int_\mathbb{R} b(t)^2\,dt
 \le C\int_D(w^2+w_\vartheta^2)
 \le C_\mu\|w\|_\mu^2.
\end{equation}
The trace is nonnegative, even, and decreasing for \(t>0\).
For completeness, angular slices converge to \(b\) in \(L^2(\mathbb{R})\) because
\(w\in H^1((a,L);L^2(\mathbb{R}))\). Each slice has the axial order, and that order is closed in \(L^2\), so it passes to the trace.
Thus
\[
 T b(T)^2\le\int_0^T b(t)^2\,dt\le C_\mu\|w\|_\mu^2
\]
for almost every \(T>0\). Since \(w(t,\vartheta)\le b(t)\), we obtain
\begin{equation}\label{eq:tail}
 \int_{|t|>T}\int_{\mathbb{S}^{N-1}}w^r
 \le C_{\mu,r}T^{-(r-2)/2}\|w\|_\mu^r
 \quad(r>2,\ T\ge1).
\end{equation}

Take a bounded sequence \(w_j\in\mathcal{K}\) and pass to a weakly convergent subsequence in \(H^1\). Its limit \(w\) belongs to \(\mathcal{K}\).
On \((-T,T)\times(\delta,L)\), ordinary two-dimensional Rellich compactness gives strong convergence in every finite \(L^r\).
The omitted caps are uniformly small by \eqref{eq:cap}, and the omitted axial ends are uniformly small by \eqref{eq:tail}.
Apply these bounds separately to \(w_j\) and \(w\), using
\(\lvert w_j-w\rvert^r\le2^{r-1}(w_j^r+w^r)\).
A diagonal subsequence converges strongly on the whole cylinder.
\end{proof}

The lemma uses both order conditions. The constants may be chosen uniformly when \(\mu\) stays in a compact subset of \((0,\infty)\) and \(r\) stays in a compact subset of \((2,\infty)\), by interpolation between two fixed exponents.

\section{Resolvent invariance at the weak level}

A cone minimum satisfies a variational inequality. To recover the full equation, we show that \(P_\mu^{-1}\) preserves the cone.
The translation and rotation fields introduced above commute with \(P_\mu\). Their weak derivatives can therefore be tested on domains bounded by reflection mirrors, without using a singular angular coordinate at the pole.

For \(Y\in\{\partial_t,Z\}\), set
\[
 \Omega_{\partial_t}=(0,\infty)\times\mathbb{S}^{N-1},
 \qquad \Omega_Z=\mathbb{R}\times Q.
\]
Functions in \(H^1_0(\Omega_Y)\) are identified with their zero extensions to \(\mathcal{C}\).

\begin{lemma}[Sign of the weak source]\label{lem:killing-source}
Let \(f\in L^2(\mathcal{C})\) be nonnegative, have the symmetries of \(\mathcal{K}\), and have its two monotonicity orders. If \(\zeta\in H^1_0(\Omega_Y)\) is nonnegative, then
\[
 \int_\mathcal{C} f\,Y\zeta\ge0.
\]
\end{lemma}

\begin{proof}
First take a nonnegative smooth test function compactly supported in \(\Omega_Y\).
The forward flow of \(\partial_t\) increases \(t\). The forward flow of \(Z\) increases \(s=\theta_N\) as long as it remains in \(Q\), since \(Zs=\theta_1>0\).
The source is nonincreasing in the corresponding coordinate: for the angular field, write \(f(t,\theta)=F(t,s)\), where \(F(t,\cdot)\) is even and nonincreasing for \(s>0\).
On the support of the test function, the flow stays in \(\Omega_Y\) for sufficiently small positive times.
Its difference quotient of \(f\) is therefore nonpositive.
Passing to the distributional limit, and using that the flow preserves volume, gives
\[
 \langle Yf,\zeta\rangle=-\int_\mathcal{C} f\,Y\zeta\le0.
\]
These inequalities concern monotone representatives and hold almost everywhere; no pointwise differentiability of \(f\) is needed.
Since \(|Y|\le1\), the functional \(\zeta\mapsto\int f\,Y\zeta\) is continuous on \(H^1_0(\Omega_Y)\).
Approximation by nonnegative compactly supported smooth functions proves the assertion.
\end{proof}

The next difference-quotient argument gives the regularity needed to test the differentiated equation. In particular, it does not assume that the source has an \(L^2\) derivative.

\begin{lemma}[Invariance of the full inverse]\label{lem:inverse}
Let \(f\in L^2(\mathcal{C})\) be nonnegative, have the symmetries of \(\mathcal{K}\), be nonincreasing in \(t>0\), and be nondecreasing in \(0<\vartheta<L\).
Then the unique \(H^1\) solution \(v=P_\mu^{-1}f\) belongs to \(\mathcal{K}\).
\end{lemma}

\begin{proof}
Coercivity gives a unique solution.
Testing with \(v^-\) gives
\(\|v^-\|_\mu^2=-\int fv^-\le0\), so \(v\ge0\).
Uniqueness gives all the required symmetries.

Fix \(Y\in\{\partial_t,Z\}\), and let \(\Phi_h\) be its isometric flow.
Because \(P_\mu\) commutes with the flow, the quotient
\(z_h=(v\circ\Phi_h-v)/h\) satisfies
\[
 \|z_h\|_\mu^2
 =\int_\mathcal{C}\frac{f\circ\Phi_h-f}{h}\,z_h
 =\int_\mathcal{C} f\,\frac{z_h\circ\Phi_{-h}-z_h}{h}.
\]
The flow estimate and \(|Y|\le1\) give
\[
 \left\|\frac{z_h\circ\Phi_{-h}-z_h}{h}\right\|_2
 \le\|Yz_h\|_2\le\|z_h\|_\mu.
\]
Thus \(\|z_h\|_\mu\le\|f\|_2\).
The quotients converge to \(z=Yv\) in \(L^2\), while a subsequence converges weakly in \(H^1\). Consequently \(Yv\in H^1(\mathcal{C})\).
Differentiating the equation weakly, or using skew-adjointness under the isometric flow, yields
\begin{equation}\label{eq:killing-weak}
 \langle z,\varphi\rangle_\mu=-\int_\mathcal{C} f\,Y\varphi
 \qquad(\varphi\in H^1(\mathcal{C})).
\end{equation}

For \(Y=\partial_t\), the function \(z\) is odd across \(t=0\).
For \(Y=Z\), it is odd across both \(\theta_1=0\) and \(\theta_N=0\): the corresponding reflections leave \(v\) invariant and reverse \(Z\).
An odd \(H^1\) function has zero trace on each mirror.
The half-cylinder and the spherical quadrant have Lipschitz boundaries, so
\[
 \zeta=z^+\mathbf 1_{\Omega_Y}\in H^1(\mathcal{C}),
 \qquad \zeta|_{\Omega_Y}\in H^1_0(\Omega_Y).
\]
In particular, the intersections of the mirrors require no additional pointwise trace.
Equation \eqref{eq:killing-weak} and Lemma~\ref{lem:killing-source} give
\[
 \|\zeta\|_\mu^2
 =\langle z,\zeta\rangle_\mu
 =-\int_\mathcal{C} f\,Y\zeta\le0.
\]
Thus \(v_t\le0\) on \(t>0\), and \(Zv\le0\) on \(\mathbb{R}\times Q\).
By \eqref{eq:killing-order}, the latter is exactly the angular order in \(\mathcal{K}\).
\end{proof}

The proof applies to every positive mass \(\mu\). The constant coefficient \(\mu^2\) does not affect commutation with either field.

\section{The minimum solves the full equation}

Write \(r=p+1>2\). The quotient in \eqref{eq:quotient} is well defined on \(\mathcal{K}\) by Lemma~\ref{lem:embedding}. Its minimum is positive, even when \(r\) is above the ordinary \(N\)-dimensional Sobolev exponent.

\begin{proposition}\label{prop:minimizer}
The infimum \(S(\mu,p)\) is positive and attained. If \(z\in\mathcal{K}\) is a minimizer with \(\int_\mathcal{C} z^r=1\), then
\[
 P_\mu z=S(\mu,p)z^p
\]
in the full distributional sense. Thus
\(w=S(\mu,p)^{1/(p-1)}z\) solves \eqref{eq:cylinder}.
\end{proposition}

\begin{proof}
The norm estimate in Lemma~\ref{lem:embedding} gives a positive lower bound. A nonzero radial test function gives a finite upper bound.
Normalize a minimizing sequence by its \(L^r\)-mass. It is bounded in \(H^1\), and weak closedness of \(\mathcal{K}\), strong \(L^r\) compactness, and lower semicontinuity produce a minimizer \(z\) with
\[
 \int_\mathcal{C} z^r=1,\qquad \|z\|_\mu^2=S:=S(\mu,p).
\]

Because \(z\in L^{2p}(\mathcal{C})\), the source \(z^p\) belongs to \(L^2\) and has the orders in Lemma~\ref{lem:inverse}. Set
\[
 v=SP_\mu^{-1}(z^p)\in\mathcal{K}.
\]
The segment \((1-\varepsilon)z+\varepsilon v\) is admissible. The right derivative of its quotient at \(\varepsilon=0\) is nonnegative, so
\begin{equation}\label{eq:cone-first-variation}
 \langle z,v-z\rangle_\mu
 -S\int_\mathcal{C} z^p(v-z)\ge0.
\end{equation}
The equation defining \(v\) gives
\[
 \langle v,z\rangle_\mu=S,\qquad
 \|v\|_\mu^2=S\int_\mathcal{C} z^pv.
\]
Substitution into \eqref{eq:cone-first-variation} yields
\(\|v\|_\mu^2\le S\).
The two preceding identities now give directly
\[
 0\le\|v-z\|_\mu^2=\|v\|_\mu^2-2S+S\le0.
\]
Thus \(v=z\), proving the full equation.
\end{proof}

We next record the regularity and strict orders that will be used in the asymptotic calculation.

\begin{lemma}\label{lem:regularity}
Every minimizing solution \(w\) is smooth and strictly positive on \(\mathcal{C}\), is bounded with a uniform local \(C^1\) bound, and satisfies \(w_t<0\) for \(t>0\).
If it is nonradial, then \(w_\vartheta>0\) on \(\mathbb{R}\times(0,L)\).
\end{lemma}

\begin{proof}
Lemma~\ref{lem:embedding} gives \(w\in L^a\) and \(w^p\in L^a\) for every finite \(a\ge2\).
Elliptic estimates on uniformly sized cylinder charts, with \(a>N\), give boundedness and a uniform \(C^1\) bound.
The nonnegative, nonzero solution is strictly positive by the strong maximum principle. Local bootstrapping then gives smoothness for any real \(p>1\).

On \(t>0\), the nonnegative function \(h=-w_t\) solves
\[
 P_\mu h=pw^{p-1}h.
\]
The local Harnack inequality makes \(h\) either identically zero or strictly positive.
The first alternative would make \(w\) constant in \(t>0\), which contradicts \(w\ne0\) and \(w\in H^1\).
For the angular order, \(g=-Zw\ge0\) on \(\mathbb{R}\times Q\) solves
\[
 P_\mu g=pw^{p-1}g,
\]
because \(Z\) commutes with \(P_\mu\).
If \(w\) is nonradial, \(g\) is not identically zero, and the same local argument gives \(g>0\).
Equation \eqref{eq:killing-order} and axial symmetry then give
\(w_\vartheta>0\) at every interior colatitude.
\end{proof}

\section{A degree-two variation and the energy gap}\label{sec:variation}

The radial homoclinic has an explicit form for any positive mass. We use a variation directed toward the equator, so it remains in the compact cone.

Put
\[
 b=\frac{\mu(p-1)}2,\qquad
 A_{\mu,p}=\left(\frac{\mu^2(p+1)}2\right)^{1/(p-1)}.
\]
The positive radial \(H^1\) solutions are translates of
\begin{equation}\label{eq:radial-cylinder}
 W(t)=A_{\mu,p}\operatorname{sech}^{2/(p-1)}(bt).
\end{equation}
Indeed the first integral of the radial equation is
\[
 (W')^2=\mu^2W^2-\frac2{p+1}W^{p+1},
\]
and integration gives \eqref{eq:radial-cylinder}. Evenness in \(t\) fixes the center.

Define the mean-zero spherical harmonic
\[
 Y_2(\theta)=1-N\theta_N^2,\qquad
 -\Delta_{\mathbb{S}^{N-1}}Y_2=2NY_2.
\]
It increases strictly from the pole toward the equator.

\begin{lemma}[An admissible negative variation]\label{lem:variation}
Let
\[
 \phi(t)=\operatorname{sech}^{(p+1)/(p-1)}(bt),\qquad
 \psi(t,\theta)=\phi(t)Y_2(\theta).
\]
For all sufficiently small \(s\ge0\), \(W+s\psi\in\mathcal{K}\).
If \(M=\int_\mathcal{C} W^{p+1}\), its quotient has zero first derivative and second derivative
\begin{equation}\label{eq:second-variation}
 2M^{-2/(p+1)}
 \left(\mu^2+2N-\frac{\mu^2(p+1)^2}{4}\right)
 \|\psi\|_2^2
\end{equation}
at \(s=0\).
\end{lemma}

\begin{proof}
The ratios
\[
 \frac{\phi}{W}=A_{\mu,p}^{-1}\operatorname{sech}(bt),\qquad
 \frac{\phi'}{W'}=\frac{p+1}{2A_{\mu,p}}\operatorname{sech}(bt)\quad(t>0)
\]
show that positivity and axial monotonicity hold, for example, when
\[
 0\le s<\frac{2A_{\mu,p}}{(p+1)(N-1)}.
\]
Angular monotonicity follows from
\(\partial_\vartheta Y_2=2N\sin\vartheta\cos\vartheta\ge0\).
All the symmetry and integrability conditions also hold.

Since \(Y_2\) has zero mean, both first variations vanish. The quotient's second derivative is
\[
 2M^{-2/(p+1)}
 \left(\|\psi\|_\mu^2-p\int_\mathcal{C} W^{p-1}\psi^2\right).
\]
Direct differentiation gives
\[
 \left(-\frac{d^2}{dt^2}+\mu^2-pW^{p-1}\right)\phi
 =\left(\mu^2-\frac{\mu^2(p+1)^2}{4}\right)\phi.
\]
Adding the angular eigenvalue \(2N\) proves \eqref{eq:second-variation}.
These differentiations are justified for every \(p>1\): the ratio
\(\psi/W\) is bounded, so the integrands and their first two derivatives are bounded by a constant times \(W^{p+1}\) near \(s=0\).
\end{proof}

\begin{proof}[Proof of Theorem~\ref{thm:cylinder}]
Proposition~\ref{prop:minimizer} and Lemma~\ref{lem:regularity} give the positive smooth solution.
Under \eqref{eq:breaking-condition}, Lemma~\ref{lem:variation} gives an admissible perturbation with quotient strictly below the value at \(W\).
A radial minimizing solution would be \(W\), by the radial classification and evenness.
The quotient is unchanged by scalar multiplication, so no normalized radial minimizer is possible.
The strict orders follow from Lemma~\ref{lem:regularity}.
\end{proof}

The same comparison gives a strict energy inequality. For a solution of the cylindrical equation, its energy is
\[
 \mathcal J_{\mu,p}(w)
 =\frac12\|w\|_\mu^2-\frac1{p+1}\int_\mathcal{C} w^{p+1}
 =\frac{p-1}{2(p+1)}\int_\mathcal{C} w^{p+1}.
\]
For a minimizing solution,
\begin{equation}
 \mathcal J_{\mu,p}(w)
 =\frac{p-1}{2(p+1)}
 S(\mu,p)^{(p+1)/(p-1)}.
\end{equation}
Therefore condition \eqref{eq:breaking-condition} implies
\begin{equation}\label{eq:energy-gap}
 \mathcal J_{\mu,p}(w)<\mathcal J_{\mu,p}(W).
\end{equation}

\section{Two ends and a positive angular coefficient}\label{sec:ends}

We first prove the exact leading decay. This step is needed both to return to the whole space and to make the projected source integrals converge.

\begin{lemma}[The leading end coefficient]\label{lem:leading}
A minimizing solution satisfies
\[
 w(t,\theta)\le Ce^{-\mu|t|}.
\]
There is \(A>0\) such that
\begin{equation}\label{eq:leading}
 w(t,\theta)=Ae^{-\mu|t|}+o(e^{-\mu|t|})
 \quad(|t|\to\infty)
\end{equation}
uniformly in direction. Its coefficient is
\[
 A=\frac1{2\mu|\mathbb{S}^{N-1}|}
 \int_\mathcal{C} e^{\mu\tau}w(\tau,\omega)^p\,d\tau\,d\omega.
\]
\end{lemma}

\begin{proof}
The trace estimate in Lemma~\ref{lem:embedding}, now for a smooth solution, gives
\(\sup_\theta w(t,\theta)\le C|t|^{-1/2}\).
Choose \(\eta\) with \(\mu/p<\eta<\mu\).
For large \(|t|\), \(w^{p-1}\le\mu^2-\eta^2\), so
\[
 (-\partial_t^2-\Delta_{\mathbb{S}^{N-1}}+\eta^2)w\le0.
\]
Comparison on each end with a multiple of \(e^{-\eta|t|}\) gives
\(w\le C_\eta e^{-\eta|t|}\).
Thus \(f=w^p\le Ce^{-\beta|t|}\), where \(\beta=p\eta>\mu\).

Use the resolvent formula
\begin{equation}\label{eq:semigroup}
 P_\mu^{-1}f=\int_0^\infty e^{-\mu^2\sigma}
 e^{\sigma\partial_t^2}e^{\sigma\Delta_{\mathbb{S}^{N-1}}}f\,d\sigma.
\end{equation}
The spherical heat semigroup is positive and contracts \(L^\infty\).
If \(F(t)=\|f(t,\cdot)\|_\infty\), then
\[
 w(t,\theta)\le\frac1{2\mu}\int_\mathbb{R}
 e^{-\mu|t-\tau|}F(\tau)\,d\tau\le Ce^{-\mu|t|}.
\]
Here the integral of the one-dimensional heat resolvent is exactly
\(e^{-\mu|t-\tau|}/(2\mu)\).
In particular the stronger source bound is
\begin{equation}\label{eq:exact-source}
 \|f(t,\cdot)\|_\infty\le Ce^{-\mu p|t|}.
\end{equation}

Let \(\Pi_0f\) be the angular average, and set
\(\bar f(t)=|\mathbb{S}^{N-1}|^{-1}\int_{\mathbb{S}^{N-1}}f(t,\theta)\,d\theta\).
The mean of the solution is
\[
 \bar w(t)=\frac1{2\mu}\int_\mathbb{R} e^{-\mu|t-\tau|}\bar f(\tau)\,d\tau.
\]
Splitting the integral at \(\tau=t\) gives
\(\bar w(t)=A_+e^{-\mu t}+O(e^{-\mu p t})\) as \(t\to+\infty\), with the stated formula for \(A_+\).
The negative-end coefficient replaces \(e^{\mu\tau}\) by \(e^{-\mu\tau}\).
Evenness of \(f\) makes the coefficients equal, and positivity makes them positive.

For the angular remainder, the first nonzero eigenvalue in the even axial class is \(2N\).
The spherical spectral gap and heat smoothing give
\begin{equation}\label{eq:heat-gap}
 \|e^{\sigma\Delta_{\mathbb{S}^{N-1}}}(I-\Pi_0)h\|_\infty
 \le C e^{-2N\sigma}\|h\|_\infty\quad(\sigma>0).
\end{equation}
For \(\sigma\ge1\), this follows by splitting off a fixed positive heat time to map \(L^2\) into \(L^\infty\), and applying the \(L^2\) spectral gap to the remaining time. For \(0<\sigma<1\), boundedness of \(I-\Pi_0\) on \(L^\infty\) gives the estimate after increasing \(C\).

Set \(\kappa_2=\sqrt{\mu^2+2N}>\mu\).
Combining \eqref{eq:semigroup}, \eqref{eq:exact-source}, and \eqref{eq:heat-gap} bounds the remainder by
\[
 C\int_\mathbb{R} e^{-\kappa_2|t-\tau|}e^{-\mu p|\tau|}\,d\tau.
\]
It is \(O(e^{-\gamma|t|})\) for every
\(\mu<\gamma<\min\{\kappa_2,\mu p\}\).
This proves \eqref{eq:leading}.
\end{proof}

The next projection determines the first angular coefficient. Its positivity follows from monotonicity rather than from a small-bifurcation calculation.

\begin{lemma}[The degree-two source is positive]\label{lem:positive-source}
If \(w\in\mathcal{K}\) is a nonradial minimizing solution, then, for every \(t\),
\[
 f_2(t):=\frac1{\|Y_2\|_2^2}
 \int_{\mathbb{S}^{N-1}}w(t,\theta)^pY_2(\theta)\,d\theta>0.
\]
Moreover,
\[
 \|Y_2\|_2^2
 =|\mathbb{S}^{N-1}|\frac{2(N-1)}{N+2}.
\]
\end{lemma}

\begin{proof}
For fixed \(t\), both \(w(t,\theta)^p\) and \(Y_2(\theta)\) increase with the colatitude from either pole to the equator, and the first is strictly increasing by Lemma~\ref{lem:regularity}.
Since \(Y_2\) has zero mean, the covariance identity gives
\[
 \int_{\mathbb{S}^{N-1}}w^pY_2
 =\frac1{2|\mathbb{S}^{N-1}|}
 \int_{\mathbb{S}^{N-1}}\int_{\mathbb{S}^{N-1}}
 (w(\theta)^p-w(\omega)^p)(Y_2(\theta)-Y_2(\omega))\,d\theta\,d\omega.
\]
The integrand is nonnegative and positive on a set of positive measure.
The norm identity follows from
\(\mathbb E[\theta_N^2]=1/N\) and
\(\mathbb E[\theta_N^4]=3/[N(N+2)]\).
\end{proof}

\begin{proposition}[Two-term cylindrical expansion]\label{prop:two-terms}
Under \eqref{eq:breaking-condition}, the minimizing solution satisfies
\begin{equation}\label{eq:two-terms}
 w(t,\theta)=Ae^{-\mu|t|}
 +Be^{-\kappa_2|t|}Y_2(\theta)
 +o(e^{-\kappa_2|t|})
\end{equation}
uniformly in direction, where \(B>0\) and
\begin{equation}\label{eq:B}
 B=\frac1{2\kappa_2\|Y_2\|_2^2}
 \int_\mathcal{C} e^{\kappa_2\tau}w(\tau,\omega)^pY_2(\omega)\,d\tau\,d\omega.
\end{equation}
When \(\mu=m=(N-2)/2\) and \(p=p_\alpha\), with \(\alpha>2\), the remainder in \eqref{eq:two-terms} is \(O(e^{-(m+4)|t|})\).
\end{proposition}

\begin{proof}
Condition \eqref{eq:breaking-condition} gives
\[
 \mu p>\frac{\mu(p+1)}2>\kappa_2.
\]
The coefficient \(w_2(t)=\|Y_2\|_2^{-2}\int w(t,\theta)Y_2(\theta)\,d\theta\) solves
\[
 -w_2''+\kappa_2^2w_2=f_2.
\]
Therefore
\[
 w_2(t)=\frac1{2\kappa_2}
 \int_\mathbb{R} e^{-\kappa_2|t-\tau|}f_2(\tau)\,d\tau.
\]
Bound \eqref{eq:exact-source} makes \eqref{eq:B} absolutely convergent.
Splitting the integral gives
\(w_2(t)=Be^{-\kappa_2|t|}+O(e^{-\mu p|t|})\).
Evenness gives the same coefficient at the two ends, and
Lemma~\ref{lem:positive-source} gives \(B>0\).

Remove the degree-zero and degree-two projections. In the remaining invariant subspace the first angular eigenvalue is \(4(N+2)\), corresponding to degree four.
The proof of \eqref{eq:heat-gap} with this spectral gap bounds the remainder by
\[
 C\int_\mathbb{R} e^{-\kappa_4|t-\tau|}e^{-\mu p|\tau|}\,d\tau,
 \qquad \kappa_4=\sqrt{\mu^2+4(N+2)}>\kappa_2.
\]
Thus it is \(O(e^{-\gamma|t|})\) for some \(\gamma>\kappa_2\).
The mean's error is \(O(e^{-\mu p|t|})\), which proves \eqref{eq:two-terms}.

For the H\'enon parameters,
\[
 \kappa_2=m+2,\qquad
 \kappa_4=m+4,\qquad
 mp_\alpha=m+\alpha+2>m+4.
\]
The last convolution is then \(O(e^{-(m+4)|t|})\), with no resonance, and the two lower projection errors decay faster.
\end{proof}

\section{Return to the critical H\'enon equation}\label{sec:henon-proof}

The cylindrical transform is exact at the exponent in \eqref{eq:henon}.
For \(u(e^t\theta)=e^{-mt}w(t,\theta)\),
\[
 m(p_\alpha+1)=N+\alpha,\qquad
 P_mw=w^{p_\alpha}.
\]
The instability condition \eqref{eq:breaking-condition} becomes
\[
 \frac{m^2(p_\alpha+1)^2}{4}>m^2+2N=(m+2)^2,
\]
which is equivalent to \(\alpha>2\).

\begin{proof}[Proof of Theorem~\ref{thm:henon}]
Apply Theorem~\ref{thm:cylinder} with \(\mu=m\) and \(p=p_\alpha\).
The resulting profile is nonradial and has the strict orders claimed.
Evenness in \(t\) gives the Kelvin identity.

Lemma~\ref{lem:leading} gives the two uniform limits and boundedness of \(u\) near zero.
On a punctured ball, \(-\Delta u=|x|^\alpha u^{p_\alpha}\) has bounded right side.
Subtract a local Newtonian potential of that bounded source. The remaining bounded harmonic function has a removable singularity at zero.
The equation therefore extends distributionally across zero.
Local \(W^{2,a}\) regularity and then Schauder estimates give
\(u\in C^2(\mathbb{R}^N)\). Its value at zero is \(A>0\).

The transform identifies the Dirichlet energy with the cylindrical norm:
\[
 \int_{\mathbb{R}^N}|\nabla u|^2
 =\int_\mathcal{C}\bigl(w_t^2+|\nabla_{\mathbb{S}^{N-1}}w|^2+m^2w^2\bigr).
\]
The cross term from \((w_t-mw)^2\) integrates to zero.
This identity extends from compactly supported smooth cylindrical functions by density in \(H^1\), proving \(u\in D^{1,2}(\mathbb{R}^N)\).
The exponent balance gives
\(\int_{\mathbb{R}^N}|x|^\alpha u^{p_\alpha+1}=\int_\mathcal{C} w^{p_\alpha+1}\).
Testing the cylindrical equation with \(w\) proves \eqref{eq:energy-identity}.
\end{proof}

The sharper angular expansion now gives the announced geometry in the original coordinates.

\begin{proof}[Proof of Theorem~\ref{thm:geometry}]
The Emden--Fowler transform preserves the energy functional on these solutions.
Inequality \eqref{eq:energy-gap} therefore gives
\(\mathcal E_\alpha(u)<\mathcal E_\alpha(U_\alpha)\);
translation of the radial profile in \(t\) does not change that energy.

Proposition~\ref{prop:two-terms} gives
\[
 w(t,\theta)=Ae^{-m|t|}
 +Be^{-(m+2)|t|}(1-N\theta_N^2)
 +O(e^{-(m+4)|t|}),\qquad B>0.
\]
Substitution at \(t\to-\infty\) and \(t\to+\infty\) proves
\eqref{eq:origin-expansion} and \eqref{eq:infinity-expansion}.
Since \(u\in C^2\), its Taylor polynomial at zero has exactly the quadratic term in \eqref{eq:origin-expansion}. The Hessian formula and the nondegenerate saddle property follow.

Finally differentiate the Kelvin identity.
The expansion at zero gives \(u(y)=A+O(|y|^2)\) and
\(\nabla u(y)=O(|y|)\).
With \(y=x/|x|^2\), its derivative has norm \(O(|x|^{-2})\), which gives the stated gradient expansion at infinity.
\end{proof}

\section{A second family in even dimension}\label{sec:second}

Let \(N=2d\) with \(d\ge2\), and let \(G_d\) be as in the introduction. For \(\theta=(\theta',\theta'')\in\mathbb{S}^{N-1}\) put \(\varphi=\arctan(|\theta''|/|\theta'|)\in[0,\pi/2]\). A \(G_d\)-invariant function on \(\mathcal{C}\) depends only on \((t,\varphi)\) and is symmetric under \(\varphi\mapsto\pi/2-\varphi\). The surface measure is a constant multiple of \(\omega(\varphi)\,d\varphi\), with \(\omega=\sin^{d-1}\varphi\cos^{d-1}\varphi\), and
\begin{equation}
 \Delta_{\mathbb{S}^{N-1}}w=w_{\varphi\varphi}+a(\varphi)w_\varphi,\qquad a=(d-1)(\cot\varphi-\tan\varphi).
\end{equation}
The orbits of \(G_d\) are products of two spheres of dimension \(d-1\), except on the strata \(\varphi=0\) and \(\varphi=\pi/2\), where they are single spheres of dimension \(d-1\).

\begin{definition}
Let \(\mathcal{K}_d\) be the set of nonnegative \(G_d\)-invariant \(w\in H^1(\mathcal{C})\) that are even in \(t\) and satisfy \(w_t\le0\) for \(t>0\) and \(w_\varphi\ge0\) for \(0<\varphi<\pi/4\), in the sense of distributions.
\end{definition}

By the symmetry, \(w_\varphi\le0\) on \((\pi/4,\pi/2)\), so a profile in \(\mathcal{K}_d\) is largest on the middle set \(\varphi=\pi/4\), where the orbits have dimension \(N-2\). Lemma~\ref{lem:embedding} holds for \(\mathcal{K}_d\), with the same proof: one takes \(D=\mathbb{R}\times(\pi/8,\pi/4)\), controls \(0<\varphi<\pi/8\) by the average over \((\pi/8,\pi/4)\) and \(\pi/4<\varphi<\pi/2\) by the symmetry, and uses the trace \(b(t)=w(t,\pi/4)\), which dominates \(w(t,\cdot)\).

The angular order of \(\mathcal{K}_d\) is not generated by a Killing field that preserves the class, so we prove the invariance of the cone by differentiating the equation. Let \(\mathcal{F}_d\) be the set of nonnegative \(G_d\)-invariant \(f\in L^2(\mathcal{C})\), even in \(t\), such that \(f(\cdot,\varphi)\) is nonincreasing on \((0,\infty)\) and \(f(t,\cdot)\) is nondecreasing on \((0,\pi/4)\) for almost every \(\varphi\) and \(t\).

\begin{lemma}\label{lem:second-inverse}
If \(f\in\mathcal{F}_d\), then \(v=P_\mu^{-1}f\) belongs to \(\mathcal{K}_d\).
\end{lemma}

\begin{proof}
As in Lemma~\ref{lem:inverse}, \(v\ge0\), \(v\) is \(G_d\)-invariant and even in \(t\), and \(v_t\le0\) for \(t>0\); that part of the proof used only the field \(\partial_t\). Expanding in the variable dual to \(t\) and in spherical harmonics, the inequality \((\xi^2+\nu+\mu^2)^2\ge\xi^4+\nu^2+2\xi^2\nu\) together with \(\int_{\mathbb{S}^{N-1}}|\nabla^2u|^2\le\int_{\mathbb{S}^{N-1}}(\Delta u)^2\) (Bochner's formula, the Ricci curvature being nonnegative) shows \(v\in H^2(\mathcal{C})\).

Write \(v=V(t,\varphi)\) and \(g=V_\varphi\). The level sets of \(\varphi\) are products of spheres of radii \(\cos\varphi\) and \(\sin\varphi\), so \(|\nabla^2_{\mathbb{S}^{N-1}}v|^2=V_{\varphi\varphi}^2+(d-1)(\cot^2\varphi+\tan^2\varphi)V_\varphi^2\). With \(d\nu=\omega\,d\varphi\,dt\) on \(\Sigma=\mathbb{R}\times(0,\pi/4)\), this gives
\begin{equation}\label{eq:second-form}
 \int_\Sigma\bigl(g_t^2+g_\varphi^2+\csc^2\varphi\,g^2\bigr)\,d\nu<\infty .
\end{equation}
Also \(g=0\) on \(\varphi=\pi/4\) in the trace sense, because \(V\) is symmetric about \(\pi/4\).

Let \(A=-\partial_\varphi^2-a\partial_\varphi\). Then \(\partial_\varphi A=(A-a')\partial_\varphi\), with \(-a'=(d-1)(\csc^2\varphi+\sec^2\varphi)\). Taking adjoints with respect to \(\omega\,d\varphi\), for \(\zeta\in C_c^\infty(\Sigma)\) and \(\psi=-\omega^{-1}(\omega\zeta)_\varphi\) we get \(P_\mu\psi=-\omega^{-1}\bigl(\omega\widetilde P\zeta\bigr)_\varphi\), where \(\widetilde P=-\partial_t^2+A-a'+\mu^2\). Testing \(P_\mu v=f\) with \(\psi\) and integrating by parts once in \(\varphi\) gives
\begin{equation}\label{eq:second-geq}
 \int_\Sigma\bigl(g_t\zeta_t+g_\varphi\zeta_\varphi+(\mu^2-a')g\zeta\bigr)\,d\nu
 =-\int_\Sigma f\,(\omega\zeta)_\varphi\,d\varphi\,dt .
\end{equation}
For \(\zeta\ge0\) the right-hand side is nonnegative, since \(f(t,\cdot)\) is nondecreasing on \((0,\pi/4)\). Both sides are continuous in the norm of \eqref{eq:second-form}, because \(|a|\le C\csc\varphi\) and \(\sec^2\varphi\le2\) on \((0,\pi/4)\). The function \(g^-\) is a limit of nonnegative elements of \(C_c^\infty(\Sigma)\) in that norm: it vanishes on \(\varphi=\pi/4\), and a cutoff \(\chi_\varepsilon(\varphi)\) near \(\varphi=0\) costs \(\int|\chi_\varepsilon'|^2g^2\,d\nu\le C\int_{\varphi<2\varepsilon}\csc^2\varphi\,g^2\,d\nu\to0\) by \eqref{eq:second-form}. Testing \eqref{eq:second-geq} with \(g^-\) gives \(-\int_\Sigma\bigl(|(g^-)_t|^2+|(g^-)_\varphi|^2+(\mu^2-a')(g^-)^2\bigr)\,d\nu\ge0\). Hence \(g^-=0\), which is the angular order.
\end{proof}

With Lemma~\ref{lem:second-inverse} in place of Lemma~\ref{lem:inverse}, Proposition~\ref{prop:minimizer} holds for the quotient \(S_d(\mu,p)\) defined as in \eqref{eq:quotient} with \(\mathcal{K}_d\) in place of \(\mathcal{K}\). Lemma~\ref{lem:regularity} also holds, with the angular part proved by the strong maximum principle for \(g=w_\varphi\ge0\), which solves \((\widetilde P-pw^{p-1})g=0\) in \(\Sigma\).

The \(G_d\)-invariant spherical harmonics are polynomials in \(|x'|^2\) and \(|x''|^2\), symmetric under the exchange. A direct computation shows that there are none of degree two or six, apart from zero, and that those of degree four are the multiples of \(H_4\). Odd degrees do not occur, since \(-I\in G_d\). On the sphere,
\[
 Y_4:=-H_4|_{\mathbb{S}^{N-1}}=\frac{d+1}d\sin^22\varphi-1,
\]
so \(Y_4\) is nondecreasing on \((0,\pi/4)\), \(|Y_4|\le1\), and \(-\Delta_{\mathbb{S}^{N-1}}Y_4=4(N+2)Y_4\). The next eigenvalue in the invariant class is at least \(8(N+6)\). Put \(\kappa_\ell=\sqrt{\mu^2+\ell(\ell+N-2)}\).

\begin{proposition}\label{prop:second-cylinder}
Let \(\mu>0\) and \(p>1\) satisfy
\begin{equation}\label{eq:second-breaking}
 \frac{\mu^2(p+1)^2}4>\mu^2+4(N+2).
\end{equation}
Then no minimizer of \(S_d(\mu,p)\) is independent of \(\theta\). Every minimizing solution \(w\) satisfies \(w_t<0\) for \(t>0\), \(w_\varphi>0\) on \(\Sigma\), and \(\mathcal J_{\mu,p}(w)<\mathcal J_{\mu,p}(W)\). There are \(A,B>0\) and \(\gamma>\kappa_4\) such that, uniformly in \(\theta\),
\[
 w(t,\theta)=Ae^{-\mu|t|}+Be^{-\kappa_4|t|}Y_4(\theta)+O(e^{-\gamma|t|}),
\]
and one can take \(\gamma=\kappa_8\) if \(\mu p>\kappa_8\).
\end{proposition}

\begin{proof}
Lemma~\ref{lem:variation} holds with \(Y_2\) and \(2N\) replaced by \(Y_4\) and \(4(N+2)\), and with \(N-1\) replaced by \(1\) in the admissible range of \(s\). The proof of Theorem~\ref{thm:cylinder} and the energy comparison then apply verbatim. Section~\ref{sec:ends} applies with the eigenvalues \(2N\) and \(4(N+2)\) replaced by \(4(N+2)\) and \(8(N+6)\), with one change: the positivity of \(\int_{\mathbb{S}^{N-1}}w^pY_4\) follows from the covariance identity, since \(w(t,\cdot)^p\) and \(Y_4\) are both functions of \(\sin^22\varphi\) and nondecreasing in it, the first strictly. Condition \eqref{eq:second-breaking} gives \(\mu p>\mu(p+1)/2>\kappa_4\).
\end{proof}

\begin{proof}[Proof of Theorem~\ref{thm:second}]
For \(\mu=m\) and \(p=p_\alpha\), condition \eqref{eq:second-breaking} reads \((N+\alpha)/2>m+4\), that is, \(\alpha>6\). Moreover \(\kappa_4=m+4\), \(\kappa_8=m+8\), and \(mp_\alpha=m+\alpha+2>m+8\). Proposition~\ref{prop:second-cylinder} therefore gives
\[
 w(t,\theta)=Ae^{-m|t|}+Be^{-(m+4)|t|}Y_4(\theta)+O(e^{-(m+8)|t|}),\qquad B>0.
\]
The return to \(\mathbb{R}^N\), the Kelvin identity, the energy identity, and \(\mathcal E_\alpha(u_2)<\mathcal E_\alpha(U_\alpha)\) follow as in the proofs of Theorems~\ref{thm:henon} and~\ref{thm:geometry}. Property (ii) is the strict angular order. Since \(r^4Y_4(\theta)=-H_4(r\theta)\), multiplying the expansion by \(e^{-mt}\) gives (iii) with \(A_2=A\) and \(B_2=B\). Finally, the expansion at the origin has no quadratic term, so \(D^2u_2(0)=0\). The solution \(u\) of Theorem~\ref{thm:henon} has \(D^2u(0)\ne0\) by Theorem~\ref{thm:geometry}, and so does \(\lambda^{(N-2)/2}u(\lambda Rx)\) for every \(\lambda>0\) and \(R\in O(N)\).
\end{proof}

\section{Compactness under parameter variation}\label{sec:parameters}

The construction applies to each parameter directly. We also show that its minimizing solutions cannot develop a concentration defect while the parameters remain in a compact subset of the positive-mass, positive-exponent region.

We first give a compactness statement that does not assume angular order. An amplitude bound and the axial normalization suffice. The radial profile \(W_{\mu,p}\) below is the function in \eqref{eq:radial-cylinder}.

\begin{lemma}[Compactness from bounded amplitudes]\label{lem:amplitude-compact}
Suppose \((\mu_j,p_j)\to(\mu,p)\in(0,\infty)\times(1,\infty)\), and let \(w_j>0\) be smooth \(H^1(\mathcal C)\) solutions of
\[
 P_{\mu_j}w_j=w_j^{p_j},\qquad
 w_j(-t,\theta)=w_j(t,\theta),\qquad
 (w_j)_t\le0\quad(t>0).
\]
Assume that \(w_j=O(e^{-\mu_j|t|})\) uniformly in direction, with constants allowed to depend on \(j\), and that \(\sup_j\|w_j\|_\infty<\infty\).
Then a subsequence converges to a positive solution \(w\) at \((\mu,p)\), strongly in \(H^1(\mathcal C)\) and in \(C^2_{\rm loc}(\mathcal C)\). Moreover,
\begin{equation}\label{eq:amplitude-relative-compact}
 \left\|\frac{w_j}{W_{\mu_j,p_j}}
             -\frac{w}{W_{\mu,p}}\right\|_\infty\longrightarrow0.
\end{equation}
There are uniform constants \(C>0\) and \(\zeta>0\) such that
\begin{equation}\label{eq:uniform-end-expansion}
 w_j(t,\theta)=A_j e^{-\mu_j|t|}
       +O\bigl(e^{-(\mu_j+\zeta)|t|}\bigr),\qquad A_j\longrightarrow A>0.
\end{equation}
\end{lemma}

\begin{proof}
The amplitude bound and local elliptic estimates give uniform \(C^{2,\gamma}\) bounds on unit cylindrical charts, for some \(\gamma>0\). Passing to a subsequence gives \(w_j\to w\) in \(C^2_{\rm loc}\), where \(w\) is nonnegative, even, and nonincreasing in \(t>0\).
Each \(w_j\) attains its maximum on \(t=0\). At a maximum point the equation gives
\[
 \|w_j\|_\infty^{p_j-1}\ge\mu_j^2.
\]
Compactness of the sphere shows that \(w\not\equiv0\), and the strong maximum principle gives \(w>0\).

For each \(j\), the Hamiltonian
\begin{equation}
 H_j(t)=\int_{\mathbb S^{N-1}}
 \left(\frac12(w_j)_t^2-\frac12|\nabla_\theta w_j|^2
       -\frac{\mu_j^2}2w_j^2+\frac1{p_j+1}w_j^{p_j+1}\right)
\end{equation}
is constant: its derivative is the equation paired with \((w_j)_t\). The assumed decay and local elliptic estimates make \(H_j(t)\to0\) at infinity, so \(H_j\equiv0\). Local convergence passes this identity to \(w\).

Since \(w(t,\theta)\) is bounded and nonincreasing in \(t>0\), it has a pointwise limit \(v(\theta)\) as \(t\to\infty\). Uniform local estimates and compactness show that every sequence of translates \(w(T+\cdot,\cdot)\), \(T\to\infty\), converges in \(C^2\) on compact slabs to the same \(t\)-independent function \(v\). In particular the convergence of the angular profiles is uniform, \(w_t\to0\), and
\[
 -\Delta_{\mathbb S^{N-1}}v+\mu^2v=v^p.
\]
Testing this equation with \(v\) and evaluating the zero Hamiltonian at the limit give
\[
 0=-\frac12\int(|\nabla_\theta v|^2+\mu^2v^2)
       +\frac1{p+1}\int v^{p+1}
   =-\frac{p-1}{2(p+1)}\int v^{p+1}.
\]
Thus \(v=0\), and \(w(t,\cdot)\to0\) uniformly at both ends.

Choose \(\mu/p<\eta<\mu\), and then \(a\) with \(\mu<a<p\eta\). For all large \(j\), \(\eta<\mu_j<a<p_j\eta\). Given a sufficiently small \(\varepsilon>0\), choose \(T\) so that \(\|w(T,\cdot)\|_\infty<\varepsilon/2\). Local convergence and monotonicity give \(w_j(t,\theta)\le\varepsilon\) on \(|t|\ge T\), uniformly for large \(j\). Choose \(\varepsilon\) so that
\(w_j^{p_j-1}\le\mu_j^2-\eta^2\) there.
Comparison on each end with a multiple of \(e^{-\eta(|t|-T)}\) yields \(w_j\le C e^{-\eta|t|}\), with a uniform constant. Hence
\[
 f_j:=w_j^{p_j}\le C e^{-a|t|}.
\]
The heat representation \eqref{eq:semigroup} and contraction of the spherical heat semigroup now give
\begin{equation}
 w_j(t,\theta)\le\frac1{2\mu_j}
    \int_{\mathbb R}e^{-\mu_j|t-\tau|}\|f_j(\tau,\cdot)\|_\infty\,d\tau
 \le C e^{-\mu_j|t|}.
\end{equation}
Local elliptic estimates give the same bound for first derivatives. Dominated convergence therefore proves strong \(H^1\) convergence.

For completeness, the angular average has leading coefficient
\[
 A_j=\frac1{2\mu_j}\int_{\mathbb R}
           e^{\mu_j\tau}\overline f_j(\tau)\,d\tau,
 \qquad
 \overline f_j=\frac1{|\mathbb S^{N-1}|}\int_{\mathbb S^{N-1}}f_j.
\]
Dominated convergence gives \(A_j\to A>0\). On mean-zero angular functions the heat gap is \(N-1\). The corresponding resolvent estimate, as in Section~\ref{sec:ends}, has decay rate \(\sqrt{\mu_j^2+N-1}>\mu_j\). Choose \(\zeta>0\) smaller than both this gap in decay rates and \(a-\mu_j\), uniformly for large \(j\). The scalar-average formula and the mean-zero estimate give \eqref{eq:uniform-end-expansion}.
Finally,
\[
 W_{\mu_j,p_j}(t)=c_j e^{-\mu_j|t|}
       \bigl(1+O(e^{-\mu_j(p_j-1)|t|})\bigr),\qquad
 c_j=2^{2/(p_j-1)}A_{\mu_j,p_j}\longrightarrow c>0.
\]
After reducing \(\zeta\), the relative tails are therefore within \(Ce^{-\zeta|t|}\) of \(A_j/c_j\). Local convergence and convergence of these coefficients prove \eqref{eq:amplitude-relative-compact}.
\end{proof}

We now apply the lemma to the cone minima, whose amplitudes are uniformly bounded by the cone embeddings and elliptic regularity.

\begin{proposition}\label{prop:parameters}
The function \(S(\mu,p)\) is continuous on \((0,\infty)\times(1,\infty)\).
Suppose \((\mu_j,p_j)\to(\mu,p)\) there, and let \(z_j\in\mathcal{K}\) be quotient minimizers normalized by
\(\int_\mathcal{C} z_j^{p_j+1}=1\).
After passage to a subsequence,
\[
 z_j\to z\quad\hbox{strongly in }H^1(\mathcal{C}),
\]
where \(z\) is a normalized minimizer at \((\mu,p)\).
The rescaled solutions \(w_j=S(\mu_j,p_j)^{1/(p_j-1)}z_j\) converge also in \(C^2_{\rm loc}(\mathcal{C})\) and satisfy
\[
 \left\|\frac{w_j}{W_{\mu_j,p_j}}-
              \frac{w}{W_{\mu,p}}\right\|_\infty\longrightarrow0.
\]
If the limiting parameters satisfy \eqref{eq:breaking-condition}, the limiting solution is nonradial.
\end{proposition}

\begin{proof}
Choose fixed \(2<r_-<p_j+1<r_+<\infty\) for all large \(j\).
Lemma~\ref{lem:embedding} and interpolation give a uniform positive lower bound for \(S(\mu_j,p_j)\).
A fixed nonzero smooth radial test profile gives a uniform upper bound.
The normalized minimizers are therefore bounded in \(H^1\).

Pass to a weakly convergent subsequence and apply cone compactness at \(r_-\) and \(r_+\).
The powers \(z_j^{r_-}\) and \(z_j^{r_+}\) converge in \(L^1\).
The bound
\[
 z_j^{p_j+1}\le z_j^{r_-}+z_j^{r_+}
\]
and almost-everywhere convergence imply
\(\int_\mathcal{C} z_j^{p_j+1}\to\int_\mathcal{C} z^{p+1}=1\).
Here the \(L^1\) convergence of the two bounding families gives both uniform integrability and tightness on the unbounded cylinder.

Since \(\mu_j\to\mu\) and \(\|z_j\|_2\) is bounded, lower semicontinuity gives
\[
 S(\mu,p)\le\|z\|_\mu^2
 \le\liminf_j S(\mu_j,p_j).
\]
Testing \(S(\mu_j,p_j)\) with a fixed normalized minimizer at \((\mu,p)\) gives the opposite upper-limit inequality by dominated convergence between two fixed finite powers.
Thus \(S(\mu_j,p_j)\to S(\mu,p)\), the limit \(z\) is a minimizer, and convergence of the Hilbert norms proves strong \(H^1\) convergence.

The rescaling factors \(S(\mu_j,p_j)^{1/(p_j-1)}\) converge. Cone embeddings at larger fixed exponents bound the solutions and their nonlinear sources uniformly in \(L^a\) for a fixed \(a>N\). Elliptic estimates on unit charts give a uniform amplitude bound. Lemma~\ref{lem:leading} gives the individual end decay, and Lemma~\ref{lem:amplitude-compact} yields both \(C^2_{\rm loc}\) convergence and the relative convergence claimed above.
Theorem~\ref{thm:cylinder} excludes a radial limiting minimizer when the strict instability condition holds.
\end{proof}

\begin{corollary}[Uniform estimates on compact parameter intervals]\label{cor:uniform-parameters}
Fix \(N\ge3\) and a nonempty compact interval \(I\subset(2,\infty)\).
There are constants \(\eta_{N,I}>0\) and \(0<c_{N,I}\le C_{N,I}<\infty\) such that, for every \(\alpha\in I\) and every entire solution \(u\) obtained from a minimizing solution for \(S(m,p_\alpha)\),
\[
 \mathcal E_\alpha(U_\alpha)-\mathcal E_\alpha(u)\ge\eta_{N,I},
 \qquad c_{N,I}\le A,B\le C_{N,I},
\]
where \(A=u(0)\) and \(B\) are the coefficients in Theorem~\ref{thm:geometry}.
In particular, every eigenvalue of the spatial Hessian \(D^2u(0)\) has absolute value at least \(2c_{N,I}\).
\end{corollary}

\begin{proof}
Every minimizing solution at a fixed parameter has energy
\[
 \frac{p_\alpha-1}{2(p_\alpha+1)}
 S(m,p_\alpha)^{(p_\alpha+1)/(p_\alpha-1)}.
\]
The radial energy and this expression depend continuously on \(\alpha\), by the explicit radial profile and Proposition~\ref{prop:parameters}. Their difference is strictly positive for \(\alpha>2\), by \eqref{eq:energy-gap}. Its minimum on \(I\) gives \(\eta_{N,I}\).

For the coefficients, take any sequence of parameters in \(I\) and corresponding minimizing profiles. After passing to a subsequence, Proposition~\ref{prop:parameters} gives a limiting minimizing solution at a parameter in \(I\). Lemma~\ref{lem:amplitude-compact} supplies a uniform bound \(w_j\le Ce^{-m|t|}\), and local convergence. Since \(\min I>2\), their sources \(w_j^{p_{\alpha_j}}\), multiplied by either \(e^{m t}\) or \(e^{(m+2)t}\), have a common integrable majorant. Dominated convergence in the coefficient formulas of Lemma~\ref{lem:leading} and \eqref{eq:B} proves \(A_j\to A_*>0\) and \(B_j\to B_*>0\). The latter positivity follows from Lemma~\ref{lem:positive-source}, since the limiting minimizer is nonradial. This proves the upper and positive lower bounds by sequential compactness. The Hessian formula in Theorem~\ref{thm:geometry} gives the last assertion.
\end{proof}

\begin{remark}\label{rem:endpoint}
Fix \(N\ge3\), let \(\alpha_j\downarrow2\), and let \(w_j\) be minimizing solutions for \(\mu=m\) and \(p=p_{\alpha_j}\). By Proposition~\ref{prop:parameters}, a subsequence converges in \(H^1(\mathcal C)\) and locally in \(C^2\) to a minimizing solution \(w_*\) at \(\alpha=2\), and \(w_j/W_{\alpha_j}\to w_*/W_2\) uniformly. If \(w_*\) is nonradial, it gives a nonradial solution at \(\alpha=2\) in this fixed dimension. If \(w_*=W_2\), local uniqueness in Corollary~\ref{cor:branch-structure} below places every sufficiently large \(j\) on the branch of Theorem~\ref{thm:local-two} with \(k=2\). In that case its cubic energy law gives
\[
 \mathcal E_{\alpha_j}(U_{\alpha_j})-\mathcal E_{\alpha_j}(u_j)
 =\mathfrak c_{N,2}(\alpha_j-2)^3+O((\alpha_j-2)^4),
\]
with the positive coefficient given in \eqref{eq:cubic-energy-constant} below. The strict local minimum at the fixed parameter \(2\), proved in Proposition~\ref{prop:local-min} below, is compatible with this latter alternative. We do not know which alternative holds globally.
\end{remark}

\section{Relative perturbations of the positive profile}\label{sec:relative}

The local argument needs a space other than the cone, since its branch below \(2\) points toward the poles. We measure perturbations relative to the radial profile, which keeps them positive and makes the nonlinear map smooth in every dimension.

For \(\alpha>-2\), define
\[
 q_\alpha=\frac{\alpha+2}{2},\qquad
 p_\alpha=1+\frac{2q_\alpha}{m},\qquad
 A_\alpha=[m(m+q_\alpha)]^{m/(2q_\alpha)},\qquad
 W_\alpha(t)=A_\alpha\operatorname{sech}^{m/q_\alpha}(q_\alpha t).
\]
Only \(\alpha\ge0\) will be used for the entire equation. The extension to \(\alpha>-2\) provides an open parameter interval at zero for the local analysis.
Let
\[
 X=\{h\in C_b^0(\mathcal{C}):h(-t,\theta)=h(t,\theta)\}
\]
with the supremum norm, and let \(X_{\rm ax}\) be its closed
\(O(N-1)\times O(1)\)-invariant subspace.
For \(\|h\|_\infty<1/2\), put
\begin{equation}\label{eq:relative-map}
 F(\alpha,h)=h-\frac1{W_\alpha}
 P_m^{-1}\!\left[
 W_\alpha^{p_\alpha}\bigl((1+h)^{p_\alpha}-1\bigr)\right].
\end{equation}
Then \(F(\alpha,0)=0\), and \(F(\alpha,h)=0\) is equivalent to
\(P_mw=w^{p_\alpha}\) for \(w=W_\alpha(1+h)>0\).
Unlike a positive-part power on a neighborhood of zero, the scalar power in \eqref{eq:relative-map} is smooth on its whole domain.

\begin{lemma}[Compactness in the relative norm]\label{lem:relative-compact}
The nonlinear term subtracted from \(h\) in \eqref{eq:relative-map} is compact on bounded subsets of \(\|h\|_\infty<1/2\), uniformly on compact parameter intervals in \((-2,\infty)\). The map \(F\) is \(C^\infty\) there, both on \(X\) and on \(X_{\rm ax}\).
\end{lemma}

\begin{proof}
We give the tail argument, including the constant angular mode.
Fix \(a>m\). If \(f\) is continuous, even in \(t\), and
\(|f(t,\theta)|\le M e^{-a|t|}\), write
\(v=P_m^{-1}f\).
Its angular average satisfies
\[
 \overline v(t)=\frac1{2m}
 \int_\mathbb{R} e^{-m|t-\tau|}\overline f(\tau)\,d\tau
 =c(f)e^{-m|t|}+O(M e^{-a|t|}),
\]
where
\[
 c(f)=\frac1{2m}\int_\mathbb{R} e^{m\tau}\overline f(\tau)\,d\tau.
\]
The same coefficient occurs at both ends because \(f\) is even.

For the mean-zero angular component, the heat semigroup on the sphere satisfies
\[
 \|e^{s\Delta_{\mathbb{S}^{N-1}}}\Pi_\perp\|_{\infty\to\infty}
 \le C e^{-(N-1)s}\qquad(s\ge0).
\]
For small \(s\), boundedness follows from the heat contraction and the bounded mean projection; for large \(s\), it follows from the spectral gap on the compact sphere. The heat representation of \(P_m^{-1}\), together with the one-dimensional heat kernel, therefore gives
\[
 \|v(t)-\overline v(t)\|_\infty
 \le C M\int_\mathbb{R} e^{-(m+1)|t-\tau|}e^{-a|\tau|}\,d\tau.
\]
We used \(\sqrt{m^2+N-1}=m+1\).
For any \(0<\eta<\min\{a-m,1\}\), it follows that
\begin{equation}\label{eq:relative-tail}
 v(t,\theta)=c(f)e^{-m|t|}
 +O(M e^{-(m+\eta)|t|})
\end{equation}
uniformly in direction.
Also,
\[
 W_\alpha(t)=C_\alpha e^{-m|t|}
 \bigl(1+O(e^{-2q_\alpha|t|})\bigr),\qquad
 C_\alpha=2^{m/q_\alpha}A_\alpha>0.
\]
Thus \(v/W_\alpha\) tends uniformly to \(c(f)/C_\alpha\) at both ends, with a uniform exponential error on compact parameter intervals.
The coefficients \(c(f)/C_\alpha\) stay bounded. Interior \(W^{2,b}\) estimates, for finite \(b>N\), give uniform local \(C^1\) bounds for \(v\). Arzel\`a--Ascoli on finite slabs, a convergent subsequence of the end coefficients, and \eqref{eq:relative-tail} prove compactness in the relative supremum norm.

To verify smoothness in the stated Banach norm, let \(E_\gamma\) be the space of continuous even sources with norm
\[
 \|f\|_{E_\gamma}=\sup_{\mathcal{C}}e^{\gamma|t|}|f(t,\theta)|.
\]
On a compact parameter interval, choose
\(m<a<\min_\alpha mp_\alpha\).
The preceding resolvent estimate gives a bounded linear map
\[
 P_m^{-1}:E_a\longrightarrow E_m\cap H^1(\mathcal{C}).
\]
For its \(H^1\) bound, use \(E_a\hookrightarrow L^2\) and coercivity of \(P_m\).

The operators \(M_\alpha:E_m\to X\), \(M_\alpha v=v/W_\alpha\), depend smoothly on \(\alpha\) in operator norm.
Indeed, \(e^{-m|t|}/W_\alpha(t)\) is uniformly bounded, and
\[
 \partial_\alpha\log W_\alpha
 =\partial_\alpha\log A_\alpha
 +\frac{m}{2q_\alpha^2}
 \bigl[\log\cosh(q_\alpha t)
       -q_\alpha t\tanh(q_\alpha t)\bigr].
\]
The bracket is bounded. Its higher parameter derivatives are also bounded on compact positive \(q_\alpha\)-intervals, since their polynomial factors in \(t\) are accompanied by exponentially decaying factors. Hence every ratio
\((\partial_\alpha^jW_\alpha)/W_\alpha\) is uniformly bounded, and the same holds for the derivatives of \(1/W_\alpha\) as multipliers from \(E_m\) to \(X\).

The source in \eqref{eq:relative-map} is a \(C^\infty\) map into \(E_a\).
Parameter derivatives of \(W_\alpha^{p_\alpha}\) are bounded by \(C(1+|t|)^je^{-mp_\alpha|t|}\), and \((1+|t|)^je^{-(mp_\alpha-a)|t|}\) is bounded since \(a<mp_\alpha\).
Also \(1+h\in(1/2,3/2)\), so all scalar derivatives of its power are bounded locally uniformly in the parameter. Taylor's formula gives the corresponding Fr\'echet derivatives in the supremum norm.
Composition with \(P_m^{-1}\) and \(M_\alpha\) proves the claimed smoothness of \(F\).
\end{proof}

We next identify the linearized operator and its range. The finite-dimensional kernel agrees with the kernel found in \cite{GGN,JX}; the relative formulation changes the nonlinear regularity required to use it.

\begin{lemma}[Kernel and range]\label{lem:relative-linear}
Put
\[
 L_\alpha=P_m-p_\alpha W_\alpha^{p_\alpha-1},\qquad
 T_\alpha h=\frac1{W_\alpha}P_m^{-1}
 (p_\alpha W_\alpha^{p_\alpha}h).
\]
Then \(D_hF(\alpha,0)=I-T_\alpha\) is Fredholm of index zero.
At \(\alpha=2(k-1)\), \(k\in\{1,2,\ldots\}\), its kernel on \(X\) consists of
\[
 \frac{\phi_k(t)Y(\theta)}{W_\alpha(t)},\qquad
 \phi_k(t)=\operatorname{sech}^{(m+k)/k}(kt),
\]
where \(Y\) runs through the degree-\(k\) spherical harmonics.
For any corresponding \(\Phi=\phi_kY\), the range annihilator is
\begin{equation}\label{eq:range-functional}
 \ell_\Phi(g)=\int_\mathcal{C} p_\alpha W_\alpha^{p_\alpha}\Phi g.
\end{equation}
The range is exactly the intersection of the kernels of these functionals.
In particular, if \(k\) is even, the kernel at \(\alpha_k\) on \(X_{\rm ax}\) is one-dimensional, spanned by \(\phi_kZ_k/W_{\alpha_k}\), where \(Z_k(\theta)=C_k^{(N-2)/2}(\theta_N)\) is the zonal harmonic of degree \(k\).
\end{lemma}

\begin{proof}
Compactness of \(T_\alpha\) follows from Lemma~\ref{lem:relative-compact}, so \(I-T_\alpha\) has index zero.
If \(h=T_\alpha h\), then
\[
 v=W_\alpha h=P_m^{-1}(p_\alpha W_\alpha^{p_\alpha}h)
 \in H^1(\mathcal{C}),\qquad L_\alpha v=0.
\]
Conversely, an even \(H^1\) kernel function with bounded ratio \(v/W_\alpha\) gives a kernel element of \(I-T_\alpha\).

Separation into spherical harmonics reduces the kernel calculation to the P\"oschl--Teller operator. Its ground eigenvalue in degree \(l\) is
\begin{equation}\label{eq:lambda-general}
 \lambda_l(\alpha)=(m+l)^2-(m+q_\alpha)^2.
\end{equation}
The next one-dimensional bound state is the odd function \(W_\alpha'\) in degree zero. All subsequent bound states, when present, have decay exponent strictly below \(m\), so cannot yield a zero eigenvalue in any angular degree. This gives the stated kernel after evenness in \(t\) removes \(W_\alpha'\).

Since \(P_m\Phi=p_\alpha W_\alpha^{p_\alpha-1}\Phi\), self-adjointness gives
\[
 \ell_\Phi(T_\alpha h)
 =\int_\mathcal{C} (P_m\Phi)P_m^{-1}
 (p_\alpha W_\alpha^{p_\alpha}h)
 =\ell_\Phi(h).
\]
All integrals converge exponentially.
Thus the range lies in the stated intersection.
The functionals are independent: on the kernel elements
\(\widetilde\Phi/W_\alpha\), their pairing is the positive definite Gram form
\(\int p_\alpha W_\alpha^{p_\alpha-1}\Phi\widetilde\Phi\).
The index-zero property gives equality of the range and that intersection.
In \(X_{\rm ax}\), the \(O(N-1)\)-invariant harmonics of degree \(k\) are the multiples of \(Z_k\), and \(Z_k\) is invariant under \(\theta_N\mapsto-\theta_N\) exactly when \(k\) is even.
\end{proof}

\section{Nonvertical branches in every dimension}\label{sec:local-branch}

Fix an even \(k\ge2\), and write \(W=W_{\alpha_k}\), \(p_k=p_{\alpha_k}\), and \(L=L_{\alpha_k}\). Put
\[
 \phi_k(t)=\operatorname{sech}^{(m+k)/k}(kt),\qquad
 \Phi=\phi_kZ_k,\qquad
 \psi=\frac{\Phi}{W}.
\]
For \(k=2\), \(Z_2=\frac{N-2}2(N\theta_N^2-1)\) is a positive multiple of \(-Y_2\), so \(\psi\) points toward the poles, opposite to the variation used in the cone argument. The local branch crosses \(\alpha_k\) because of the sign of a cubic moment.

\begin{lemma}\label{lem:cubic}
For every even \(k\ge2\), \(\int_{\mathbb{S}^{N-1}}Z_k^3>0\). For \(k=2\),
\[
 \frac1{|\mathbb{S}^{N-1}|}\int_{\mathbb{S}^{N-1}}(N\theta_N^2-1)^3
 =\frac{8(N-1)(N-2)}{(N+2)(N+4)},\qquad
 \frac1{|\mathbb{S}^{N-1}|}\int_{\mathbb{S}^{N-1}}(N\theta_N^2-1)^2
 =\frac{2(N-1)}{N+2}.
\]
\end{lemma}

\begin{proof}
Let \(\lambda=(N-2)/2>0\). Rogers' linearization formula \cite{AAR} writes \((C_k^\lambda)^2\) as a combination of \(C_{2k-2j}^\lambda\), \(0\le j\le k\), with positive coefficients. Since \(k\) is even, \(C_k^\lambda\) occurs, for \(j=k/2\). Zonal harmonics of different degrees are orthogonal on the sphere, so \(\int Z_k^3\) is a positive multiple of \(\int Z_k^2\). The values for \(k=2\) follow from the averages \(1/N\), \(3/(N(N+2))\), and \(15/(N(N+2)(N+4))\) of \(\theta_N^2\), \(\theta_N^4\), and \(\theta_N^6\).
\end{proof}

\begin{proof}[Proof of Theorem~\ref{thm:local-two}]
At \((\alpha_k,0)\), Lemma~\ref{lem:relative-linear} gives a simple kernel \(\operatorname{span}\{\psi\}\) on \(X_{\rm ax}\). Using
\[
 (I-T_\alpha)h=\frac1{W_\alpha}P_m^{-1}L_\alpha(W_\alpha h),
\]
and applying \(\ell_\Phi\), we obtain
\begin{equation}\label{eq:relative-transversality}
 \ell_\Phi\bigl(\partial_\alpha D_hF(\alpha_k,0)\psi\bigr)
 =\int_\mathcal{C}\Phi(\partial_\alpha L_\alpha|_{\alpha=\alpha_k})\Phi
 =-(m+k)\int_\mathcal{C}\Phi^2<0.
\end{equation}
Indeed, the derivative of \(W_\alpha h\) contributes a term annihilated by \(L\Phi=0\), and differentiation of \eqref{eq:lambda-general} at \(q_\alpha=k\) gives the last equality. The Crandall--Rabinowitz theorem \cite{CR} gives a smooth curve
\[
 \alpha(0)=\alpha_k,\qquad h(s)=s\psi+O(s^2)\quad\hbox{in }X_{\rm ax}.
\]
Its profiles are positive because \(\|h(s)\|_\infty<1/2\).

Differentiating \(F(\alpha(s),h(s))=0\) twice and applying \(\ell_\Phi\) gives
\begin{equation}\label{eq:alpha-slope}
 \alpha'(0)
 =-\frac{p_k(p_k-1)}{2(m+k)}
 \frac{\displaystyle\int_\mathcal{C} W^{p_k-2}\Phi^3}
 {\displaystyle\int_\mathcal{C}\Phi^2}.
\end{equation}
The numerator is \(\int_\mathbb{R}W^{p_k-2}\phi_k^3\,dt\cdot\int_{\mathbb{S}^{N-1}}Z_k^3\). The first factor converges, since \(W^{p_k-2}\phi_k^3\) decays like \(e^{-(5k+2m)|t|}\), also when \(p_k<2\), and the second is positive by Lemma~\ref{lem:cubic}. Consequently \(\alpha'(0)<0\). For \(k=2\) the ratio is explicit: with \(I(a)=\int_\mathbb{R}\operatorname{sech}^a(2t)\,dt=\tfrac12\mathrm B(\tfrac12,\tfrac a2)\),
\[
 \frac{\int W_2^{p_2-2}\Phi^3}{\int\Phi^2}
 =\frac{N-2}2\,A_2^{p_2-2}\frac{I(m+5)}{I(m+2)}\frac{4(N-2)}{N+4}.
\]
Local inversion of \(s\mapsto\alpha(s)\) covers a punctured interval on both sides of \(\alpha_k\). The degree-\(k\) projection of \(h(s)\) is nonzero for small \(s\ne0\), so the profiles are nonradial.

Every profile satisfies \(w=P_m^{-1}(w^{p_\alpha})\), with an \(L^2\) source, and therefore belongs to \(H^1(\mathcal{C})\). It is smooth and bounded by \(Ce^{-m|t|}\). The scalar-average formula used in Lemma~\ref{lem:relative-compact}, now with the positive source \(w^{p_\alpha}\), gives positive leading coefficients at both ends, equal by evenness. The return to the entire equation, removability at zero, finite energy, and the Kelvin identity follow as in the proof of Theorem~\ref{thm:henon}.

It remains to compare energies. Put
\[
 J=\int_\mathcal{C} W^{p_k-2}\Phi^3>0,\qquad M=\int_\mathcal{C}\Phi^2.
\]
Along the solution curve, \(v_s=w_s-W_{\alpha(s)}\) satisfies
\[
 v_s=s\Phi+R_s,\qquad \|R_s\|_{H^1}\le Cs^2,\qquad |R_s|\le Cs^2W_{\alpha(s)}.
\]
The pointwise estimate follows from the expansion in \(X_{\rm ax}\) and the bound \(|W_\alpha/W_{\alpha_k}-1|\le C|\alpha-\alpha_k|\). The \(H^1\) estimate follows by differentiating the resolvent equation into \(H^1\), using its \(E_a\)-valued source and the bounded map in Lemma~\ref{lem:relative-compact}.

Write \(\alpha=\alpha(s)\) and \(p=p_\alpha\). Relative positivity permits Taylor expansion of the scalar power even when \(p<2\):
\[
 \mathcal E_\alpha(u_s)-\mathcal E_\alpha(U_\alpha)
 =\frac12\langle L_\alpha v_s,v_s\rangle
 -\frac{p(p-1)}6\int_\mathcal{C} W_\alpha^{p-2}v_s^3
 +\mathcal R_s,
\]
where the first pairing is the quadratic form on \(H^1\), and
\[
 |\mathcal R_s|\le C\int_\mathcal{C} W_\alpha^{p+1}\Bigl|\frac{v_s}{W_\alpha}\Bigr|^4\le Cs^4.
\]
The constants are uniform near \(\alpha_k\). Put \(\tau=\alpha(s)-\alpha_k=O(s)\). Because \(L\Phi=0\) and the potential depends smoothly on \(\alpha\) in \(L^\infty\), \(\|L_\alpha\Phi\|_2\le C|\tau|\). Thus the cross term involving \(s\Phi\) and \(R_s\), and the quadratic term in \(R_s\), are \(O(s^4)\). The eigenvalue identity gives \(\langle L_\alpha\Phi,\Phi\rangle=\lambda_k(\alpha)M+O(\tau^2)\), and the exponential bounds of Lemma~\ref{lem:relative-compact} give \(\int W_\alpha^{p-2}v_s^3=s^3J+O(s^4)\). Consequently
\[
 \mathcal E_{\alpha(s)}(u_s)-\mathcal E_{\alpha(s)}(U_{\alpha(s)})
 =\frac12s^2\lambda_k(\alpha(s))M-\frac{p_k(p_k-1)}6Js^3+O(s^4).
\]
Since \(\lambda_k(\alpha(s))=-(m+k)\alpha'(0)s+O(s^2)\), \eqref{eq:alpha-slope} gives
\[
 \mathcal E_{\alpha(s)}(u_s)-\mathcal E_{\alpha(s)}(U_{\alpha(s)})
 =\frac{p_k(p_k-1)}{12}Js^3+O(s^4).
\]
For \(s>0\) we have \(\alpha(s)<\alpha_k\) and a positive energy difference. For \(s<0\) both signs reverse.
\end{proof}

The part of the branch with \(\alpha<\alpha_k\) has higher energy than the radial profile. For \(k=2\) it is therefore not a cone minimizer, which has energy at most \(\mathcal E_\alpha(U_\alpha)\). The argument gives no lower bound for \(\varepsilon_{N,k}\).

The slope also identifies the complete local solution set in this symmetry class and determines the small eigenvalue of the linearization. We record these consequences with the energy written in the original parameter.

\begin{corollary}[Local uniqueness, energy, and Morse index]\label{cor:branch-structure}
Fix \(N\ge3\) and an even \(k\ge2\). In a neighborhood of \((\alpha_k,0)\) in \(\mathbb R\times X_{\rm ax}\), all solutions of \(F(\alpha,h)=0\) belong either to the radial family or to the branch in Theorem~\ref{thm:local-two}. In particular, for each sufficiently close \(\alpha\ne\alpha_k\) there is exactly one nonradial profile in that neighborhood.

Put
\[
 M_k=\int_{\mathcal C}\Phi^2,\qquad
 J_k=\int_{\mathcal C}W_{\alpha_k}^{p_k-2}\Phi^3>0,
 \qquad \Phi=\phi_kZ_k.
\]
Along this branch,
\begin{equation}\label{eq:cubic-energy-law}
 \mathcal E_\alpha(u_\alpha)-\mathcal E_\alpha(U_\alpha)
 =-\mathfrak c_{N,k}(\alpha-\alpha_k)^3
       +O(|\alpha-\alpha_k|^4),
\end{equation}
where
\begin{equation}\label{eq:cubic-energy-constant}
 \mathfrak c_{N,k}
 =\frac{2(m+k)^3M_k^3}{3p_k^2(p_k-1)^2J_k^2}>0.
\end{equation}
The linearization \(P_m-p_\alpha w_\alpha^{p_\alpha-1}\) is nondegenerate on the \(t\)-even, \(O(N-1)\times O(1)\)-invariant subspace. Its Morse index on that subspace is
\begin{equation}
 \operatorname{ind}_{\rm ax,e}(w_\alpha)=
 \begin{cases}
 k/2+1,&\alpha<\alpha_k,\\
 k/2,&\alpha>\alpha_k.
 \end{cases}
\end{equation}
The index counts negative \(L^2\) eigenvalues with multiplicity in this subspace; it is not the unrestricted Morse index.
\end{corollary}

\begin{proof}
The local completeness assertion is part of the simple-eigenvalue bifurcation theorem \cite{CR}. Since \(\alpha'(0)\ne0\), its curve is a graph over \(\alpha\). Substitution of
\[
 s=-\frac{2(m+k)M_k}{p_k(p_k-1)J_k}
           (\alpha-\alpha_k)+O((\alpha-\alpha_k)^2)
\]
into the energy expansion just proved gives \eqref{eq:cubic-energy-law}--\eqref{eq:cubic-energy-constant}. The coefficient is unchanged if \(\Phi\) is multiplied by any nonzero constant.

At \(\alpha=\alpha_k\), separation of variables shows that the even axial linearization has exactly \(k/2\) negative eigenvalues: one ground state for each angular degree \(0,2,\ldots,k-2\). Its kernel is simple, spanned by \(\Phi\). The next radial bound state is odd in \(t\); all subsequent ones are positive, as in Lemma~\ref{lem:relative-linear}. Thus the rest of the spectrum is separated from zero.

The potentials \(p_{\alpha(s)}w_s^{p_{\alpha(s)}-1}\) depend smoothly on \(s\) in \(L^\infty\), also when \(p_k<2\), by relative positivity and exponential decay. They vanish uniformly at the ends. Bounded self-adjoint perturbation theory therefore gives a unique simple eigenvalue \(\rho(s)\) near zero, with \(\rho(0)=0\); all other signs remain fixed. Differentiating at zero and testing against \(\Phi\) gives
\[
 \rho'(0)=-(m+k)\alpha'(0)-\frac{p_k(p_k-1)J_k}{M_k}
          =(m+k)\alpha'(0).
\]
Consequently
\[
 \rho(\alpha)=(m+k)(\alpha-\alpha_k)
                    +O((\alpha-\alpha_k)^2).
\]
It is negative below \(\alpha_k\) and positive above. This proves the index formula and nondegeneracy in the stated subspace.
\end{proof}

The coefficient at \(k=2\) controls the energy gap of every global cone minimizer, even if such minimizers do not approach the radial profile at the endpoint.

\begin{corollary}[A cubic lower bound for the cone energy gap]\label{cor:cubic-gap}
Let \(u_\alpha\) be any minimizing solution given by Theorem~\ref{thm:henon}. For \(\alpha>2\) sufficiently close to \(2\),
\begin{equation}\label{eq:cubic-cone-gap}
 \mathcal E_\alpha(U_\alpha)-\mathcal E_\alpha(u_\alpha)
 \ge \mathfrak c_{N,2}(\alpha-2)^3-C_N(\alpha-2)^4,
\end{equation}
where \(C_N<\infty\). Equality holds up to \(O((\alpha-2)^4)\) if the minimizers converge to \(W_2\); see Remark~\ref{rem:endpoint}.
\end{corollary}

\begin{proof}
Write \(\tau=\alpha-2>0\), \(\kappa=m+2\), and
\[
 \widetilde\Phi_\alpha(t,\theta)
 =-\operatorname{sech}^{(p_\alpha+1)/(p_\alpha-1)}(q_\alpha t)Z_2(\theta).
\]
This points toward the equator. By Lemma~\ref{lem:variation},
\(z_{\alpha,a}=W_\alpha+a\widetilde\Phi_\alpha\) belongs to \(\mathcal K\) for all small \(a\ge0\), uniformly for \(\alpha\) near \(2\). Put \(D=p_2(p_2-1)J_2>0\), where \(p_2=p_{\alpha=2}\). Taylor expansion in the relative norm and the exact ground eigenvalue yield
\[
 \mathcal J_{m,p_\alpha}(z_{\alpha,a})
       -\mathcal J_{m,p_\alpha}(W_\alpha)
 =-\frac{\kappa M_2}2\tau a^2+\frac D6a^3
       +O(\tau^2a^2+\tau a^3+a^4).
\]
Here \(\widetilde\Phi_2=-\Phi\), so its cubic moment is \(-J_2\), and
\(\lambda_2(\alpha)=-\kappa\tau-\tau^2/4\).

Let \(c_{\alpha,a}>0\) maximize \(\mathcal J_{m,p_\alpha}(cz_{\alpha,a})\) over \(c>0\). Both linear mass variations vanish because \(Z_2\) has zero angular mean, and hence \(c_{\alpha,a}=1+O(a^2)\). The maximization changes the preceding expansion by \(O(a^4)\). Taking \(a=2\kappa M_2\tau/D\) gives
\[
 \mathcal J_{m,p_\alpha}(c_{\alpha,a}z_{\alpha,a})
       -\mathcal J_{m,p_\alpha}(W_\alpha)
 =-\mathfrak c_{N,2}\tau^3+O(\tau^4).
\]
For any nonzero \(z\in\mathcal K\), its maximum energy along the positive scalar ray is
\[
 \max_{c>0}\mathcal J_{m,p}(cz)
 =\frac{p-1}{2(p+1)}
 \left(\frac{\|z\|_m^2}{(\int z^{p+1})^{2/(p+1)}}\right)^{(p+1)/(p-1)}.
\]
A cone minimizing solution has no greater energy than this test value. The Emden--Fowler energy identity proves \eqref{eq:cubic-cone-gap}.
\end{proof}

At \(\alpha=2\) the second variation \eqref{eq:second-variation} vanishes in the cone direction, and the cubic moment of Lemma~\ref{lem:cubic} decides the local picture inside the cone.

\begin{proposition}\label{prop:local-min}
Let \(\alpha=2\), \(p=p_2\), and put
\[
 R(w)=\frac{\|w\|_m^2}{\left(\int_{\mathcal C}w^{p+1}\right)^{2/(p+1)}}
 \qquad(0\ne w\in\mathcal K).
\] There is \(\rho>0\) such that, if \(w\in\mathcal{K}\), \(c>0\), and \(\|w/(cW_2)-1\|_\infty<\rho\), then \(R(w)\ge R(W_2)\), with equality only if \(w\) is a multiple of \(W_2\).
\end{proposition}

\begin{proof}
Write \(W=W_2\), \(M_0=\int W^{p+1}\), \(L=P_m-pW^{p-1}\), and \(\Phi=\phi_2Y_2\), with \(Y_2=1-N\theta_N^2\) as in Section~\ref{sec:variation}. Since \(R\) is invariant under scaling, we may replace \(c\) by \(\int W^pw/M_0\). For small \(\rho\) this keeps \(|h|<3\rho\), where \(v=w/c-W\) and \(h=v/W\), and it gives \(\int W^pv=0\), that is, \(\langle W,v\rangle_m=0\).

Write \(v=a\Phi+v_\perp\) with \(\int W^{p-1}\Phi v_\perp=0\), and put \(n=\|v_\perp\|_m\). For each \(t\), \(v(t,\cdot)\) is nondecreasing in \(\vartheta\) on \((0,\pi/2)\) and even under the reflection, while \(Y_2\) has mean zero and is strictly increasing there. The covariance identity gives \(\int_{\mathbb{S}^{N-1}}v(t,\theta)Y_2(\theta)\,d\theta\ge0\), hence \(a\ge0\). Also \(a\le C\rho\), and \(|v_\perp|\le C\rho W\).

The quadratic form \(Q(z)=\langle Lz,z\rangle\) satisfies \(Q(z)\ge c_0\|z\|_m^2\) on invariant functions, even in \(t\), with \(\langle W,z\rangle_m=\langle\Phi,z\rangle_m=0\). Indeed, on functions of \(t\) alone \(W\) minimizes the one-dimensional quotient, so \(Q\ge0\) on \(\{\langle W,z\rangle_m=0\}\), and its translation kernel is spanned by the odd function \(W'\). A zero of the constrained radial form satisfies \(Lz=cW^p\); since \(LW=-(p-1)W^p\), adding a multiple of \(W\) reduces it to this kernel, and the tangent condition forces that multiple to vanish. Evenness then excludes \(W'\). In degree two the kernel is spanned by \(\Phi\), and in degree \(\ell\ge4\) the lowest eigenvalue is \((m+\ell)^2-(m+2)^2>0\). In each sector \(L\) is a compact perturbation of \(P_m\), which gives coercivity. Since \(L\Phi=0\) and \(p\int W^{p-1}\Phi v_\perp=\langle\Phi,v_\perp\rangle_m=0\), we get \(Q(v)=Q(v_\perp)\ge c_0n^2\).

Put \(S_j=\int W^{p+1}h^j\). Since \(\|W+v\|_m^2=M_0+\|v\|_m^2\) and \(|h|<1/2\), Taylor expansion of \(\int W^{p+1}(1+h)^{p+1}\) and of the power \(-2/(p+1)\) gives
\[
 \frac{R(w)}{R(W)}-1=\frac1{M_0}\Bigl(Q(v)-\frac{p(p-1)}3S_3\Bigr)+E,
 \qquad
 |E|\le C\bigl(S_2^2+S_2\|v\|_m^2+\textstyle\int W^{p+1}h^4\bigr).
\]
Now \(S_2\le C(a^2+\min\{n^2,\rho^2\})\), \(\|v\|_m^2\le C(a^2+n^2)\), and \(\int W^{p+1}h^4\le C(a^4+\rho^2n^2)\), so \(|E|\le C(a^4+(a+\rho)n^2)\). In \(S_3\) the leading term is \(a^3\int W^{p-2}\Phi^3\), which is negative by Lemma~\ref{lem:cubic}, since \(Y_2\) is a negative multiple of \(Z_2\). The other terms are bounded by \(Ca^2n+C(a+\rho)n^2\), because \(W^{p-2}\Phi^2\in L^2(\mathcal{C})\), \(\Phi/W\) is bounded, and \(|v_\perp|\le C\rho W\). Hence, with \(c_3=-\frac{p(p-1)}3\int W^{p-2}\Phi^3>0\),
\[
 M_0\Bigl(\frac{R(w)}{R(W)}-1\Bigr)\ge\frac{c_0}2n^2-Ca^2n+c_3a^3-Ca^4
 \ge\frac{c_0}4n^2+\frac{c_3}2a^3
\]
for small \(\rho\). The right-hand side vanishes only if \(a=n=0\), that is, \(v=0\).
\end{proof}

Proposition~\ref{prop:local-min} does not decide which alternative in Remark~\ref{rem:endpoint} holds. Relative compactness now identifies a radial limiting family with the local branch, but strict local minimality at the fixed parameter \(2\) does not exclude its lower-energy members for \(\alpha>2\). The quadratic coefficient becomes negative there, and Corollary~\ref{cor:cubic-gap} quantifies the resulting energy gain.

\section{Local obstructions at the two endpoints}\label{sec:local-endpoints}

At \(\alpha=2\) the obstruction in odd dimension is algebraic. At \(\alpha=0\) the radial solution is surrounded by the translated bubbles, and this manifold has to be factored out before the implicit function theorem applies.

\subsection{Dilation normalization}

Evenness in the cylindrical variable can always be achieved by a dilation.

\begin{lemma}\label{lem:t-symmetry}
Let \(w>0\) solve \(P_mw=w^p\), with \(p>1\), be smooth, and satisfy
\(w=O(e^{-m|t|})\) uniformly in direction. Then, for some \(t_0\), it is even about \(t=t_0\) and strictly decreasing for \(t>t_0\).
\end{lemma}

\begin{proof}
Write \(w_\lambda(t,\theta)=w(2\lambda-t,\theta)\) for \(t<\lambda\).
On the set where \(w>w_\lambda\), the coefficient
\((w^p-w_\lambda^p)/(w-w_\lambda)\) is at most \(p w^{p-1}\).
For \(\lambda\) sufficiently negative, this is less than \(m^2/2\).
A negative-part test for \(w_\lambda-w\) on the half-cylinder therefore starts the comparison. Exponential decay and local elliptic estimates justify its \(H^1\) cutoffs.

Slide the plane to its limiting value \(\lambda_*\).
It is finite, since comparison for every \(\lambda\) would make \(w\) nondecreasing in \(t\), contradicting its two zero limits.
At \(\lambda_*\), the strong maximum principle gives either equality or strict comparison.
In the latter case, choose \(R>0\) with \(-R<\lambda_*-1\) so that
\(p w^{p-1}<m^2/2\) for \(t<-R\), and choose a small \(\varepsilon>0\).
Strict comparison on the compact slab
\([-R,\lambda_*-\varepsilon]\times\mathbb{S}^{N-1}\) persists for nearby planes.
Hence, for \(\lambda>\lambda_*\) sufficiently close, any possible negative part is supported in the separated regions
\(t<-R\) and \(\lambda_*-\varepsilon<t<\lambda\).
It has zero trace at the intervening slab and at the reflection plane.
On the tail, its energy test gives coercivity with mass at least \(m^2/2\).
On the strip, if \(C=\sup_\mathcal{C} pw^{p-1}\), the same test gives
\[
 \int (|\nabla z|^2+m^2z^2)\le C\int z^2
 \le C\frac{(\lambda-\lambda_*+\varepsilon)^2}{\pi^2}
       \int |z_t|^2 .
\]
Choosing \(\varepsilon\) and then \(\lambda-\lambda_*\) small forces both negative parts to vanish. This contradicts the definition of \(\lambda_*\).
Thus \(w=w_{\lambda_*}\).
The preceding comparisons give \(w_t\le0\) for \(t>\lambda_*\).
The function \(-w_t\) solves
\(P_m(-w_t)=pw^{p-1}(-w_t)\) there.
Local Harnack makes it strictly positive, since it cannot vanish identically for a nonzero profile with zero limit at infinity.
\end{proof}

A translation of \(t\) is a dilation of \(u\). In particular, every positive entire solution with the stated Newtonian decay has such an even normalization.

\subsection{The fixed value 2 in odd dimension}

At \(\alpha=2\), all degree-two harmonics belong to the kernel on \(X\).
Represent them by
\[
 Y_A(\theta)=\theta^{\mathsf T}A\theta,\qquad
 A=A^{\mathsf T},\qquad \operatorname{tr}A=0.
\]
The relevant quadratic map is
\[
 Q(A)=A^2-\frac{\operatorname{tr}(A^2)}N I.
\]

\begin{proof}[Proof of Proposition~\ref{prop:local-rigidity}(i)]
Let \(h\) be a sufficiently small solution of \(F(2,h)=0\).
Split
\[
 h=\frac{\phi Y_A}{W_2}+\eta
\]
using a bounded projection onto the finite-dimensional kernel.
The index-zero range description in Lemma~\ref{lem:relative-linear} makes the linearized map an isomorphism from a fixed complementary domain onto its range.
Projecting the equation onto this range and using its bounded inverse yields
\[
 \|\eta\|_\infty\le C\|h\|_\infty^2.
\]
Hence \(\|h\|_\infty\le C|A|\) and \(\|\eta\|_\infty\le C|A|^2\) after reducing the neighborhood; if \(A=0\), the same estimate already gives \(h=0\).
Here \(|A|^2=\operatorname{tr}(A^2)\).

For \(v=W_2h\), the equation gives
\[
 L_2v=\frac{p_2(p_2-1)}2 W_2^{p_2-2}v^2
 +O(W_2^{p_2}|h|^3).
\]
Pair this identity with \(\phi Y_B\), for an arbitrary symmetric trace-free \(B\). The left side vanishes. The preceding estimates and exponential decay give
\[
 0=\frac{p_2(p_2-1)}2
 \left(\int_\mathbb{R} W_2^{p_2-2}\phi^3\right)
 \int_{\mathbb{S}^{N-1}}Y_A^2Y_B+O(|A|^3|B|).
\]
The sixth spherical moment is
\[
 \int_{\mathbb{S}^{N-1}}Y_A^2Y_B
 =\frac{8|\mathbb{S}^{N-1}|}{N(N+2)(N+4)}
 \operatorname{tr}(A^2B).
\]
Equivalently,
\[
 \Pi_2(Y_A^2)=\frac4{N+4}Y_{Q(A)}.
\]
Since the radial factor is positive, we conclude
\[
 |Q(A)|\le C|A|^3.
\]

If \(N\) is odd, \(Q(A)=0\) forces \(A=0\).
Indeed, for nonzero \(A\) it would imply \(A^2=cI\), with \(c>0\); its eigenvalues would all be \(\pm\sqrt c\), and their sum could not vanish in odd dimension.
By compactness of the unit sphere in the trace-free symmetric matrices,
\[
 |Q(A)|\ge c_N|A|^2
\]
for some \(c_N>0\).
The two bounds force \(A=0\) when \(|A|\) is small, proving \(h=0\).
\end{proof}

Even dimension has an actual family at the fixed parameter, already obtained in \cite[Proposition 1.7]{GGN}. If \(N=2d\), put
\(H(x)=|x'|^2-|x''|^2\), with \(x',x''\in\mathbb{R}^d\).
In the present coefficient convention and Kelvin normalization, the family is
\begin{equation}\label{eq:even-explicit}
 u_b(x)=[(N+2)(N-2)]^{(N-2)/8}
 (1-b^2)^{(N-2)/8}
 [1+|x|^4-2bH(x)]^{-(N-2)/4},
 \quad |b|<1.
\end{equation}
It is positive and nonradial for \(b\ne0\).
Its energy equals that of \(U_2\), since differentiating the energy along this fixed-exponent solution family gives zero.
Thus the odd-dimensional local statement does not extend to even dimension.
As \(b\to1\), \(\max u_b\to\infty\) and \(u_b\to0\) locally uniformly away from the sphere \(\{|x'|=1,\ x''=0\}\). In the companion paper \cite{LeConc} the same kind of sphere appears as a concentration orbit near the parameter \(2(k-1)/(N-k-1)\), for \(O(k)\times O(N-k)\)-invariant solutions with \(k<N/2\). When \(k=N/2\), this parameter equals \(2\), and the algebraic factor \(2k-N\) in the formal geometric balance vanishes. The exact family \eqref{eq:even-explicit} describes concentration at this fixed endpoint parameter.

\subsection{A uniform radial neighborhood near zero}

At \(\alpha=0\), the even cylindrical kernel is the degree-one space.
Let
\[
 \phi_1(t)=\operatorname{sech}^{m+1}t,\qquad
 \psi_a(t,\theta)=\frac{\phi_1(t)(a\cdot\theta)}{W_0(t)}.
\]
Translation invariance produces exact solutions near \(W_0\).
For \(|b|<1\), these are
\begin{equation}\label{eq:yamabe-family}
 w_b(t,\theta)=A_0
 \left(\frac{\sqrt{1-|b|^2}}{\cosh t-b\cdot\theta}\right)^m.
\end{equation}
They are even in \(t\), solve \(P_mw_b=w_b^{p_0}\), and satisfy
\[
 \frac{w_b}{W_0}-1
 =m\operatorname{sech} t\,(b\cdot\theta)+O(|b|^2)
 \quad\hbox{in }X.
\]
The derivative spans the full kernel.

\begin{proof}[Proof of Proposition~\ref{prop:local-rigidity}(ii)]
Apply finite-dimensional Lyapunov--Schmidt reduction to \(F\) near \((0,0)\), using Lemmas~\ref{lem:relative-compact} and \ref{lem:relative-linear}.
Write
\[
 h=\psi_a+\eta(\alpha,a),\qquad
 \eta(\alpha,0)=0,\qquad
 D_a\eta(0,0)=0,
\]
where \(\eta\) lies in a fixed complement of the kernel.
The projected range equation determines a smooth \(\eta\).
For \(\Phi_j=\phi_1\theta_j\), let
\(\ell_j=\ell_{\Phi_j}\) in \eqref{eq:range-functional}, evaluated at \(\alpha=0\).
These are fixed functionals, independent of the nearby parameter.
Define the remaining equation by
\[
 G_j(\alpha,a)=\ell_j
 \bigl(F(\alpha,\psi_a+\eta(\alpha,a))\bigr),\quad 1\le j\le N.
\]
The range description makes \(G=0\) equivalent to the full equation once the range equation has been solved. Also \(G(\alpha,0)=0\).

The exact family \eqref{eq:yamabe-family} implies
\(G(0,a)=0\) for every small \(a\).
To justify this for the chosen kernel coordinate, project
\(w_b/W_0-1\) onto the kernel. Its derivative at zero is invertible, so these projections cover a full neighborhood of \(a=0\); uniqueness in the range equation identifies their complementary parts with \(\eta(0,a)\).
Hadamard's lemma now gives
\[
 G(\alpha,a)=\alpha\,\widehat G(\alpha,a)
\]
with \(\widehat G\) smooth.

The derivative \(D_a\widehat G(0,0)\) is invertible.
Put \(\mathcal A_\alpha=D_hF(\alpha,0)\).
Since \(\ell_i\mathcal A_0=0\) and \(D_a\eta(0,0)=0\), differentiation gives
\[
 (D_a\widehat G(0,0))_{ij}
 =\ell_i(\partial_\alpha\mathcal A_\alpha|_0\,\psi_{e_j}).
\]
Terms containing \(\partial_\alpha D_a\eta(0,0)\) are annihilated by \(\ell_i\mathcal A_0\). As in \eqref{eq:relative-transversality}, the resulting matrix is
\[
 \left[
 \int_\mathcal{C}(\phi_1\theta_i)
 (\partial_\alpha L_\alpha|_{\alpha=0})
 (\phi_1\theta_j)
 \right]_{i,j=1}^N
 =-(m+1)\left(\int_\mathbb{R}\phi_1^2\right)
 \frac{|\mathbb{S}^{N-1}|}{N}I.
\]
Here \(\lambda_1'(0)=-(m+1)\) in \eqref{eq:lambda-general}.
The implicit-function theorem for \(\widehat G\) therefore gives the unique small solution \(a=0\) for all small \(\alpha\).
Since \(\widehat G(\alpha,0)=0\), this uniqueness is uniform across \(\alpha=0\).
For \(0<\alpha<\delta_N\), the original reduced equation is equivalent to \(\widehat G=0\), which proves the assertion with a fixed relative radius \(\rho_N\).
\end{proof}

Thus nonradial solutions for small \(\alpha>0\) cannot be obtained as small relative perturbations of the radial one. Solutions far from it are not excluded.

\section{An application with a Hardy potential}

Critical equations with an inverse-square potential have an established bifurcation theory. Dancer, Gladiali, and Grossi \cite{DGG} studied the unweighted critical equation, and Bonheure, Casteras, and Gladiali \cite{BCG} treated singular weights \(|x|^{-s}\), \(0\le s<2\). The following application concerns positive H\'enon weights. The free positive mass in the cylindrical construction gives a solution directly at each parameter satisfying the cone instability condition, without changing the cone. Let
\[
 m=\frac{N-2}{2},\qquad \lambda<m^2,\qquad
 \mu=\sqrt{m^2-\lambda}.
\]
The Emden--Fowler transform takes the equation
\begin{equation}\label{eq:hardy}
 -\Delta u-\frac{\lambda}{|x|^2}u
 =|x|^\alpha u^{p_\alpha}
\end{equation}
to \(P_\mu w=w^{p_\alpha}\).

\begin{corollary}
Let \(N\ge3\), \(\alpha>0\), and
\[
 \lambda<\Lambda_N(\alpha)
 :=m^2\left(1-\frac{8N}{(\alpha+2)(\alpha+2N-2)}\right).
\]
Then \eqref{eq:hardy} has a positive nonradial weak solution in
\(D^{1,2}(\mathbb{R}^N)\), smooth on \(\mathbb{R}^N\setminus\{0\}\).
For some \(A>0\), uniformly in direction,
\[
 \lim_{r\to0}r^{m-\mu}u(r\theta)=A,\qquad
 \lim_{r\to\infty}r^{m+\mu}u(r\theta)=A.
\]
Its weighted nonlinear mass and Hardy integral are finite.
\end{corollary}

\begin{proof}
The parameter condition is exactly
\[
 \mu^2(p_\alpha-1)(p_\alpha+3)>8N,
\]
which is equivalent to \eqref{eq:breaking-condition}.
Theorem~\ref{thm:cylinder} and Lemma~\ref{lem:leading} apply.

The transform gives \(u\in D^{1,2}\),
\[
 \int_{\mathbb{R}^N}\frac{u^2}{|x|^2}=\int_\mathcal{C} w^2<\infty,\qquad
 \int_{\mathbb{R}^N}|x|^\alpha u^{p_\alpha+1}
 =\int_\mathcal{C} w^{p_\alpha+1}<\infty.
\]
It solves the equation on the punctured space.
Near zero, \(u=O(r^{-m+\mu})\).
Both \(|x|^{-2}u\) and \(|x|^\alpha u^{p_\alpha}\) are locally integrable, since their radial integral exponents are respectively \(m+\mu-1\) and \(m+p_\alpha\mu-1\).
Testing with a cutoff around zero and letting its radius tend to zero extends the weak equation across zero: the derivative-cutoff error is bounded by
\(\|\nabla u\|_{L^2(B_{2\varepsilon})}O(\varepsilon^{(N-2)/2})\).
The asymptotic formulas and positivity follow from the cylindrical result.
\end{proof}

For \(\lambda=0\), the condition in this corollary is exactly \(\alpha>2\).
For
\(\lambda<-m^2(N+1)/(N-1)\), it holds for every \(\alpha>0\). This concerns a repulsive potential and says nothing about \eqref{eq:henon} for \(\alpha\le2\).

\section{The lower-parameter obstruction}

For every \(\alpha>0\) the radial solution has a negative direction of degree one, but it is odd under the equatorial reflection and does not belong to \(\mathcal{K}\). A balance identity shows that no solution can be monotone from one pole to the other.

\begin{proposition}[Balance and one-pole rigidity]
Every positive classical solution of \eqref{eq:henon} with Newtonian decay and \(\alpha>0\) satisfies
\begin{equation}\label{eq:translation-balance}
 \int_{\mathbb{R}^N}|x|^{\alpha-2}x\,u^{p_\alpha+1}\,dx=0.
\end{equation}
If, for some axis, \(u(r\theta)\) is axially symmetric and nonincreasing in colatitude on the whole interval \(0<\vartheta<\pi\) for every \(r>0\), then \(u\) is radial.
\end{proposition}

\begin{proof}
Interior estimates on rescaled annuli give
\(\nabla u=O(|x|^{1-N})\) at infinity.
Multiply the equation by \(\partial_i u\) on
\(B_R\setminus B_\varepsilon\) and integrate by parts.
The two terms from the Dirichlet energy are boundary terms and tend to zero.
The nonlinear term becomes
\[
 \int |x|^\alpha u^{p_\alpha}\partial_i u
 =-\frac{\alpha}{p_\alpha+1}
 \int |x|^{\alpha-2}x_i u^{p_\alpha+1}+o(1).
\]
Its boundary terms also vanish. The displayed vector integral is absolutely convergent at zero and infinity, by boundedness near zero and Newtonian decay. Letting \(R\to\infty\) and \(\varepsilon\to0\), and using \(\alpha>0\), proves \eqref{eq:translation-balance}.

In the stated angular class, both \(\theta_N=\cos\vartheta\) and
\(u(r\theta)^{p_\alpha+1}\) are nonincreasing in \(\vartheta\).
The covariance identity for the spherical angular measure gives
\[
 \int_{\mathbb{S}^{N-1}}\theta_N u(r\theta)^{p_\alpha+1}\,d\theta
 \ge0,
\]
with strict inequality whenever the angular profile is nonconstant.
The \(N\)-th component of \eqref{eq:translation-balance} is its integral against the positive radial measure \(r^{N+\alpha-2}\,dr\).
Thus every angular profile is constant, first for almost every radius and then for every radius by continuity.
\end{proof}

The same monotone class also defeats minimization, for a different reason.

\begin{proposition}
Let \(\mathcal{K}_{\rm pole}\) have all the conditions in Definition~\ref{def:cone}, except that \(w_\vartheta\le0\).
For \(r>2N/(N-2)\),
\[
 \inf_{\substack{w\in\mathcal{K}_{\rm pole}\cap L^r(\mathcal{C})\\w\ne0}}
 \frac{\|w\|_\mu^2}{\left(\int_\mathcal{C} w^r\right)^{2/r}}=0.
\]
The same conclusion holds in the one-pole cone: nonnegative \(O(N-1)\)-invariant \(H^1\) profiles, even and nonincreasing in \(t>0\), with \(w_\vartheta\le0\) on the full interval \(0<\vartheta<\pi\), without equatorial reflection.
\end{proposition}

\begin{proof}
Choose a nonzero smooth nonnegative nonincreasing cutoff \(\chi\) on \([0,\infty)\), supported in \([0,1)\), and set
\[
 z_\varepsilon(t,\vartheta)
 =\varepsilon^{-N/r}
 \chi\!\left(\frac{t^2+\vartheta^2}{\varepsilon^2}\right).
\]
This function is supported near the north pole and belongs to the one-pole cone.
Adding its reflection across the equator gives a profile in \(\mathcal{K}_{\rm pole}\).
Using \(\sin\vartheta\sim\vartheta\) and changing variables gives
\[
 \int_\mathcal{C} z_\varepsilon^r\to c_1>0,\qquad
 \|z_\varepsilon\|_\mu^2
 =c_2\varepsilon^{N-2-2N/r}(1+o(1))
 +O(\varepsilon^{N-2N/r}),
\]
with \(c_2>0\).
The first exponent is positive, so the quotient tends to zero.
Reflection merely multiplies the two integrals by fixed positive factors and gives the same conclusion for \(\mathcal{K}_{\rm pole}\).
\end{proof}

At the critical H\'enon exponent,
\[
 N-2-\frac{2N}{p_\alpha+1}
 =\frac{(N-2)\alpha}{N+\alpha}>0.
\]
Thus this obstruction occurs for every \(\alpha>0\).

The rigidity near zero can be strengthened: a nonradial family cannot stay bounded in the relative norm as \(\alpha\to0\).

\begin{proposition}[Escape from bounded relative sets]\label{prop:escape}
Suppose \(\alpha_j\downarrow0\), and \(u_j\) are positive nonradial entire solutions with Newtonian decay. Normalize their cylindrical profiles \(w_j\) to be even by dilation. Then
\[
 \left\|\frac{w_j}{W_{\alpha_j}}\right\|_\infty\longrightarrow\infty.
\]
\end{proposition}

\begin{proof}
Suppose a subsequence has a uniform bound \(M\) in that norm, and put
\(R_j=w_j/W_{\alpha_j}\).
The compactness proof in Lemma~\ref{lem:relative-compact} also applies to the sources
\(W_{\alpha_j}^{p_{\alpha_j}}R_j^{p_{\alpha_j}}\) for
\(0\le R_j\le M\); smoothness near zero is not needed for this compactness assertion.
Since
\[
 R_j=\frac1{W_{\alpha_j}}P_m^{-1}
 (W_{\alpha_j}^{p_{\alpha_j}}R_j^{p_{\alpha_j}}),
\]
we may assume \(R_j\to R_*\) uniformly.
The resolvent bound gives
\(\|R_j\|_\infty\le C\|R_j\|_\infty^{p_{\alpha_j}}\), with \(C\) uniform.
Thus \(R_*\not\equiv0\).
The limit \(w_*=W_0R_*\) is a positive even \(H^1\) solution at \(\alpha=0\).
After return to \(\mathbb{R}^N\), removability and the classification theorem of Caffarelli--Gidas--Spruck \cite{CGS} identify it as a Yamabe bubble.
Its even normalization makes it one of the profiles \(w_b\) in \eqref{eq:yamabe-family}.

In cylindrical coordinates, \eqref{eq:translation-balance} reads
\[
 \int_\mathcal{C} e^{-t}\theta\,w_j^{p_{\alpha_j}+1}\,dt\,d\theta=0,
\]
because \(m(p_\alpha+1)=N+\alpha\).
The uniform relative bound gives an integrable majorant at both ends, independent of \(j\), so this identity passes to \(w_b\).
If \(b\ne0\), the profile \(w_b(t,\theta)^{p_0+1}\) is strictly increasing in \(b\cdot\theta\). Covariance on the sphere then gives
\[
 \int_\mathcal{C} e^{-t}(b\cdot\theta)w_b^{p_0+1}>0,
\]
a contradiction.
Hence \(b=0\), \(w_*=W_0\), and \(R_j\to1\) uniformly.
Proposition~\ref{prop:local-rigidity}(ii) now forces \(w_j=W_{\alpha_j}\) for all large \(j\), contradicting nonradiality.
\end{proof}

The preceding relative escape combines with bounded-amplitude compactness to give the dilation-invariant conclusion stated in the introduction.

\begin{proof}[Proof of Theorem~\ref{thm:amplitude}]
If \eqref{eq:amplitude-blowup} fails, pass to a subsequence with
\(\mathcal A(u_j)\le M\). Lemma~\ref{lem:t-symmetry} allows a dilation for which every \(w_j\) is even and decreasing in \(t>0\). The dilation leaves \(\mathcal A(u_j)\) unchanged.
Classical regularity at zero and Newtonian decay give \(w_j=O(e^{-m|t|})\) at both ends. Local elliptic estimates give the corresponding gradient decay, hence \(w_j\in H^1(\mathcal C)\).
Apply Lemma~\ref{lem:amplitude-compact} with \(\mu_j=m\) and \(p_j=p_{\alpha_j}\to p_0\). Its relative convergence makes \(w_j/W_{\alpha_j}\) uniformly bounded, contradicting Proposition~\ref{prop:escape}.

For a fixed \(M\), failure of the first assertion would give positive nonradial solutions with \(0<\alpha_j<1/j\) and \(\mathcal A(u_j)\le M\). This contradicts the conclusion just proved.
\end{proof}

An additional consequence concerns the variational construction itself. Although the cone contains nonradial functions, its minimum stays radial in a whole interval next to zero.

\begin{corollary}[Radial cone minima near zero]\label{cor:small-cone}
For each \(N\ge3\), there is \(\delta_N>0\) such that, for \(0<\alpha<\delta_N\), every minimizing solution for \(S(m,p_\alpha)\) equals \(W_\alpha\). Equivalently, for every \(v\in\mathcal K\),
\[
 \|v\|_m^2\ge
 \left(\int_{\mathcal C}W_\alpha^{p_\alpha+1}\right)^{(p_\alpha-1)/(p_\alpha+1)}
 \left(\int_{\mathcal C}v^{p_\alpha+1}\right)^{2/(p_\alpha+1)},
\]
with equality for nonzero \(v\) exactly when \(v=cW_\alpha\), \(c>0\).
\end{corollary}

\begin{proof}
If nonradial minimizing solutions existed for a sequence \(\alpha_j\downarrow0\), the bounds in Proposition~\ref{prop:parameters} would give uniformly bounded cylindrical amplitudes. The return argument in Section~\ref{sec:henon-proof}, which uses only \(\alpha>0\) for removability and finite energy, gives positive classical entire solutions with Newtonian decay. This contradicts Theorem~\ref{thm:amplitude}. Hence all minimizing solutions are radial for sufficiently small positive \(\alpha\), and the radial classification identifies them with \(W_\alpha\). The sharp inequality and its equality case follow from the quotient definition and \(\|W_\alpha\|_m^2=\int W_\alpha^{p_\alpha+1}\).
\end{proof}

Taken together, these results give nonradial solutions for every \(\alpha>2\), for \(\alpha\) in a dimension-dependent interval below \(2\), and at \(\alpha=2\) in even dimension. The concentration constructions developed in \cite{LeConc} concern additional lower-parameter intervals and are not used in any proof above. The union of the available local constructions has not been shown to cover \((0,2)\), and unrestricted existence at \(\alpha=2\) in odd dimension remains unresolved. The following appendix gives an algebraic complement to the local rigidity statement at that fixed parameter.

\appendix
\section{Polynomial denominators at the parameter two}

The local rigidity result in Proposition~\ref{prop:local-rigidity}(i) concerns arbitrary sufficiently small relative perturbations. There is also a global algebraic statement within a particular class: the reciprocal of a polynomial of degree at most four. The Newtonian condition fixes its leading homogeneous part, and the equation excludes all cubic and linear terms. This calculation does not impose an amplitude restriction, but its polynomial assumption is essential.

\begin{proposition}[The complete polynomial-denominator test]
Suppose \(D\) is a positive real polynomial of degree at most four and
\[
 u=C D^{-(N-2)/4},\qquad C>0,
\]
is a positive solution of \eqref{eq:henon} at \(\alpha=2\), with positive uniform Newtonian coefficient at infinity.
After multiplication of \(D\) by a positive constant and absorption of this constant in \(C\), necessarily
\[
 D=c+|x|^4-2x^{\mathsf T}Ax,\qquad
 A=A^{\mathsf T},\quad\operatorname{tr}A=0,
\]
\[
 A^2=\gamma I,\qquad
 C^{p_2-1}=(N^2-4)(c-\gamma),\qquad 0\le\gamma<c.
\]
Conversely these data give a solution.
In odd dimension every such solution is radial.
In even dimension, the nonradial members are precisely the balanced eigenvalue family of \cite[Proposition 1.7]{GGN}, up to rotation and dilation.
\end{proposition}

\begin{proof}
Put \(b=(N-2)/4\), so \(bp_2=b+2\).
Newtonian decay implies that \(D\) has degree four and its leading part is a positive multiple of \(|x|^4\); normalize this multiple to one.
Write \(D=|x|^4+P_3+P_2+P_1+c\) in homogeneous parts.
Clearing denominators gives the polynomial identity
\begin{equation}\label{eq:polynomial-denominator-identity}
 bD\Delta D-b(b+1)|\nabla D|^2=C^{p_2-1}|x|^2.
\end{equation}
Its degree-six terms cancel because \(4(N+2)=16(b+1)\).
Its degree-five term is
\[
 b|x|^2\bigl(|x|^2\Delta P_3-2(N+2)P_3\bigr).
\]
Decompose \(P_3=H_3+|x|^2H_1\), with harmonic homogeneous parts.
Since \(\Delta(|x|^2H_1)=2(N+2)H_1\), the displayed term forces \(H_3=0\).
Thus \(P_3=|x|^2\ell\cdot x\).
Let \(h=\ell/4\), and translate \(y=x+h\).
The translated denominator has no cubic part and can be written
\[
 \widetilde D=|y|^4-2y^{\mathsf T}Ay+d\cdot y+c'.
\]
The right side of \eqref{eq:polynomial-denominator-identity} is now
\(C^{p_2-1}|y-h|^2\).
Its left side has degree-four term
\(-4b\operatorname{tr}A\,|y|^4\), which forces \(\operatorname{tr}A=0\).
Its remaining degree-three term is
\(8b(b+1)|y|^2d\cdot y\), which forces \(d=0\).
The left side is then the homogeneous quadratic
\[
 (N^2-4)y^{\mathsf T}(c'I-A^2)y.
\]
The right side has constant term \(C^{p_2-1}|h|^2\).
Hence \(h=0\), and the original cubic part is zero.
Comparing the remaining quadratic forms gives the conditions in the statement.
Trace freeness and \(A^2=\gamma I\) imply that the eigenvalues are \(\pm\sqrt\gamma\), with equal multiplicities unless \(\gamma=0\).
Odd dimension therefore forces \(A=0\).
Positivity of the minimum denominator gives \(c>\gamma\), and direct substitution proves the converse.
\end{proof}

This exhausts the polynomial-denominator ansatz with the required decay. It does not prove radiality of all positive solutions at \(\alpha=2\) in odd dimension.

\section*{Statements and Declarations}
%
%\noindent\textbf{Funding.} This research received no external funding.
%
%\smallskip
\noindent\textbf{Competing interests.} The author declares no conflict of interest.

\smallskip
\noindent\textbf{Data availability.} No data were generated or analyzed in this study.

%\smallskip
%\noindent\textbf{Use of AI tools.} AI-assisted tools were used during the preparation of the manuscript. The author reviewed and edited all AI-assisted output and takes full responsibility for the content of the work.

\end{document}